\documentclass[smallextended]{svjour3}

\usepackage{listings}
\usepackage{fix-cm}
\usepackage{graphicx}
\usepackage[a4paper, total={7in, 10in}]{geometry}
\usepackage{multirow}
\usepackage{amsmath,amsfonts,color}
 
\usepackage{comment}
\usepackage{hyperref}
\usepackage{authblk}
\hypersetup{colorlinks=true,linkcolor=blue,citecolor=blue,filecolor=black,urlcolor=MidnightBlue}    
\usepackage{mathtools}
\newtheorem{thm}{Theorem}

\begin{document}

\title{An Efficient Quasi-Newton Method with Tensor Product Implementation for Solving Quasi-Linear Elliptic Equations and Systems}

\author{Wenrui Hao         \and
        Sun Lee \and
        Xiangxiong Zhang
}


\institute{Wenrui Hao \at
                Department of Mathematics, The Pennsylvania State University, University Park, PA 16802, USA \\
              \email{wxh64@psu.edu}    \and
              Sun Lee \at
                Department of Mathematics, The Pennsylvania State University, University Park, PA 16802, USA \\
              \email{skl5876@psu.edu} 
           \and
              Xiangxiong Zhang \at
                Department of Mathematics, Purdue University, West Lafayette, IN 47907, USA \\
              \email{zhan1966@purdue.edu}
}

\maketitle

\begin{abstract}
    In this paper, we introduce a quasi-Newton method optimized for efficiently solving quasi-linear elliptic equations and systems, with a specific focus on GPU-based computation. By approximating the Jacobian matrix with a combination of linear Laplacian and simplified nonlinear terms, our method reduces the computational overhead typical of traditional Newton methods while handling the large, sparse matrices generated from discretized PDEs. We also provide a convergence analysis demonstrating local convergence to the exact solution under optimal choices for the regularization parameter, ensuring stability and efficiency in each iteration. Numerical experiments in two- and three-dimensional domains validate the proposed method’s robustness and computational gains with tensor product implementation. This approach offers a promising pathway for accelerating quasi-linear elliptic equations and systems solvers, expanding the feasibility of complex simulations in physics, engineering, and other fields leveraging advanced hardware capabilities.

\end{abstract}

\section{Introduction}

Numerical solutions of partial differential equations (PDEs) are central to modeling and simulating physical phenomena across a broad range of fields, including physics, engineering, and finance. The process of solving PDEs numerically involves discretizing the continuous domain and operators, transforming them into systems of algebraic equations that can be addressed computationally \cite{leveque2007finite,smith1985numerical}. Common discretization techniques, such as finite difference methods (FDM), finite element methods (FEM), and finite volume methods (FVM) \cite{leveque2007finite,perrone1975general,strikwerda2004finite}, yield nonlinear systems that approximate the PDE solutions on a discrete set of points.

For quasi-linear elliptic PDEs, solving these nonlinear systems requires significant computational resources, particularly when iterative methods are employed to resolve large, sparse linear systems \cite{hao2020homotopy,hao2024companion,wang2018two}. Newton's method and its variants are widely used to address these systems due to their local quadratic convergence properties \cite{kelley2003solving,lu2020non}, but their computational cost can be prohibitive, especially for high-dimensional or complex nonlinear systems, due to the need for repeated linearizations and inversion of large Jacobian matrices \cite{dembo1982inexact,kelley1987quasi,liu2019sparsifying,saad2003iterative}. Quasi-Newton methods, which approximate the Jacobian matrix, have shown promise for reducing computational expense while maintaining convergence rates, particularly in high-dimensional or multi-physics applications \cite{kelley1987quasi,li2023approximate,lu2020full,nocedal2006large}. 

The computational burden of Newton-type methods is compounded in large-scale problems, where memory and processing requirements are substantial. Recent advances in computational hardware, especially the development of graphics processing units (GPUs), provide new opportunities for accelerating numerical computations using parallel processing on tensor structures \cite{liu2024simple,nickolls2010gpu,owens2008gpu}. Despite this, the need to repeatedly compute and invert large nonlinear Jacobians in Newton’s method presents challenges for effective GPU parallelization.

In this work, we propose an efficient quasi-Newton method tailored for quasi-linear elliptic equations that leverages GPU architectures and tensor product structures to enhance computational efficiency. Our approach approximates the Jacobian matrix by utilizing the linear Laplacian operator along with simplified representations of the nonlinear terms, effectively transforming the quasi-linear problem into a form more amenable to parallel computation and GPU acceleration.

The quasi-Newton method developed here is motivated by the observation that many discretization schemes lead to linear systems characterized by sparse matrices, where nonlinear terms can often be approximated by diagonal or block-diagonal matrices \cite{liu2019sparsifying,saad2003iterative}. Discretizing the equations using appropriate numerical methods leads to a nonlinear system
\[
    A_h U_h + N_h(U_h) = 0,
\]
where \( A_h \) is the discretized Laplacian operator or linear operator, and \( N_h \) is the discretized nonlinear terms with numerical solution \(U_h\). Traditional Newton's methods at each iteration a linear system involving the Jacobian matrix \( J_h = A_h + N_u(U_h) \), where \( N_u \) is the Jacobian of the nonlinear term. Newton's method emphasize direct inversion of these operators \(J_h\), which becomes computationally demanding as matrix size and complexity increase. Our quasi-Newton approach incorporates a scaled identity matrix as a proxy for the nonlinear terms in the Jacobian
\[
    [A_h+\beta I]^{-1} \approx [A_h+N_u]^{-1},
\]
thereby simplifying the inversion process and enabling more efficient computation. This approach also ensures convergence with proper assumptions, as discussed in \S \ref{IM}.

The primary reason for the reduction in computational effort is that we only need to compute the diagonalization of matrix \( A_h \) once, and then we can reuse them repeatedly. If we can express \( A_h \) in its diagonalized form as \( A_h = T \Lambda T^{-1} \), where \( T \) contains the eigenvectors and \( \Lambda \) is the diagonal matrix of eigenvalues, then it follows that
\[
A_h + \beta I = T (\Lambda + \beta ) T^{-1}.
\]
To obtain the inverse of \( A_h + \beta I \), we need only compute the inverse of the diagonal matrix \( \Lambda + \beta I \), which is simple since it involves inverting the diagonal entries. Moreover, during the quasi-Newton process, we only update \( \beta \), so at each iterations, we simply compute the inverse of the diagonal matrix \(\Lambda + \beta I\). 

An additional reason for the efficiency of our approach lies in the ability to exploit the structure of tensor products and Kronecker products. This allows us to compute the inverse of matrices like \( T (\Lambda + \beta I) T^{-1} \) even more efficiently by leveraging properties of these products to avoid full matrix inversion. Specifically, the Kronecker product structure lets us perform on smaller matrices, significantly reducing the computational, memory cost and making the method more amenable to GPU acceleration. Specifically, when \( A_h \) arises from discretizations on tensor-product grids, it can often be expressed as a sum of Kronecker products of smaller matrices. For example, in two dimensions, the matrix \( A_h \) can be written as
\begin{equation*}
    A_h = I_y \otimes A_x + A_y \otimes I_x,
\end{equation*}
where \( A_x \) and \( A_y \) are one-dimensional discretization matrices, and \( I_x \), \( I_y \) are identity matrices of appropriate sizes with the Kronecker product \( \otimes \). By diagonalizing the smaller matrices \( A_x \) and \( A_y \), say
\[
A_x = T_x \Lambda_x T_x^{-1}, \quad A_y = T_y \Lambda_y T_y^{-1},
\]
we can express \( A_h \) using the properties of Kronecker products:
\begin{align*}
    A_h + \beta I &= (I_y \otimes A_x) + (A_y \otimes I_x) + \beta (I_y \otimes I_x) \\
    &= (T_y \otimes T_x) \left( I_y \otimes \Lambda_x + \Lambda_y \otimes I_x + \beta I_y \otimes I_x \right) (T_y^{-1} \otimes T_x^{-1}).
\end{align*}
The matrix inside the parentheses is diagonal matrix, and its inverse involves only element-wise inverse:
\begin{equation}\label{eq:inverse_tensor}
    [A_h + \beta I]^{-1} = (T_y \otimes T_x) \left( I_y \otimes (\Lambda_x + \beta I) + (\Lambda_y + \beta I) \otimes I_x \right)^{-1} (T_y^{-1} \otimes T_x^{-1}).
\end{equation}
Since \( \Lambda_x \) and \( \Lambda_y \) are diagonal matrices, the inverse in \eqref{eq:inverse_tensor} reduces to element-wise inversion, which can be performed efficiently. This exploitation of the tensor and Kronecker product structures significantly reduces computational cost and is highly amenable to parallelization on GPUs.

In this paper, our primary objective is to harness the computational power of modern GPUs using a minimalist \textsc{Matlab} implementation, originally developed by \cite{liu2024simple}. The code from \cite{liu2024simple} is both simple and efficient, making it an ideal foundation for our work. To adapt this code to our problem, we reformulated our setting using tensor products. The main challenge arose due to the presence of a nonlinear term, which precluded the direct application of the existing method. However, our new method, quasi-Newton method, overcomes this obstacle by approximating the nonlinear term with a simple identity matrix scaled by a constant multiplier.

The proposed quasi-Newton method is compatible with various discretization techniques and problem dimensions. In this paper, we apply the Spectral Element Method (SEM) to two- and three-dimensional problems. While both SEM and the Finite Difference Method (FDM) are well-suited for GPU acceleration, SEM is emphasized here due to its greater accuracy and complexity.
\textcolor{black}{On a rectangular mesh for a rectangular domain,  the $Q^k$ element SEM is obtained by implementing classical Lagrangian $Q^k$ basis continuous finite element method (FEM)  with $(k+1)$-point Gauss-Lobatto quadrature. This particular choice of quadrature not only results in a diagonal mass matrix but also makes SEM become a finite difference type scheme defined at quadrature points, allowing easier implementations, especially on GPUs.
When regarded as a finite difference scheme, 
the {\it a priori} error estimate of  $Q^k$ element SEM ($k\geq 2$)
is $(k+2)$-th order accurate in $\ell^2$-norm  for second order PDEs
\cite{li2022accuracy,li2020superconvergence}. In other words, the standard superconvergence result of functions values at Gauss-Lobatto points of classical FEM still hold in this simple implementation. Moreover,
SEM is  more robust than the standard FEM in the sense of provable monotonicity. Recently, $Q^2$ element SEM has been proven monotone for an elliptic operator with scalar coefficient \cite{cross2024-Q2,li2020monotonicity}, and $Q^3$ element SEM can be proven monotone for Laplacian in two dimensions \cite{cross2023monotonicity}, with applications to quite a few second order PDEs \cite{hu2023positivity,liu2024structure,liu2023positivity,shen2022discrete}.  See \cite{cross2023monotonicity} for numerical results about monotonicity of $Q^k$ element SEM ($k\geq 4$) for Laplacian.   
}

We organize the paper as follows: In \S \ref{Model}, we introduce the problem setup and basic assumptions. In \S \ref{Discret_2d}, we describe the spectral element method in two and three dimensions. We introduce the quasi-Newton method in \S \ref{IM} and discuss its implementation. Convergence proofs under specific assumptions are provided in \S \ref{Convergence}. Finally, in \S \ref{examples}, we demonstrate the performance of our method on several numerical examples. 

\section{Model Problem}\label{Model}

We consider the following quasi-linear elliptic equation defined in a bounded rectangular domain $\Omega \subset \mathbb{R}^d$:

\begin{equation}\label{eq:main-pde}
  -\Delta u + \lambda u + f(x, u) = 0, \quad \text{in } \Omega \hbox{~and~}
\color{black}{\alpha_1 \frac{\partial u}{\partial n}+\alpha_2 u=\alpha_3}, \quad \text{on } \partial\Omega.
\end{equation}

Here, $\Delta$ denotes the Laplacian operator, $\lambda \in \mathbb{R}$ is a constant parameter, $u: \Omega \rightarrow \mathbb{R}$ is the unknown function to be determined, and $f: \Omega \times \mathbb{R} \rightarrow \mathbb{R}$ is a given nonlinear function \textcolor{black}{with $f(x,~\cdot~)\in C^2(\mathbb{R})$}. The symbol $\partial u / \partial \boldsymbol{n}$ represents the normal derivative of $u$ on the boundary $\partial \Omega$, with $\boldsymbol{n}$ being the outward unit normal vector. We have Dirichlet boundary conditions if \textcolor{black}{$\alpha_1=0$ and  Neumann boundary conditions if $\alpha_2=0$.}

Eq.~\eqref{eq:main-pde} models a variety of physical phenomena, including steady-state heat conduction, electrostatics, and diffusion processes with reaction terms \cite{kapitula2004counting,kevrekidis2015defocusing,simon1995concentration}. The presence of the nonlinear function $f(x, u)$ introduces additional complexity, necessitating advanced numerical methods for its solution.


\section{Numerical Discretization}\label{Discret_2d}
In this section, we describe two discretization methods with tensor product structures that result in a more efficient computational implementation. In this case, the proposed Quasi-Newton Method with this particular structure shows efficiency with GPU implementation.

\subsection{Spectral Element Method in 2D}\label{subsec:spectral-element-2d}

We consider the spectral element method on a rectangular mesh with continuous piecewise \( Q^k \) polynomial basis functions in \cite{liu2024simple}. This involves the tensor product of piecewise polynomials of degree \( k \), utilizing \((k + 1)\)-point Gauss-Lobatto quadrature on rectangular meshes. In a two-dimensional element, the basis functions are products of one-dimensional piecewise polynomials of degree \( k \), and all integrals are approximated using Gauss points. For a rectangular mesh using the \( Q^k \) basis, let \( (x_i, y_j) \) denote all the nodal points, where \( i = 1, \dotsc, N_x \) and \( j = 1, \dotsc, N_y \).

\begin{figure}[ht]
    \centering
    \includegraphics[width=.2\textwidth]{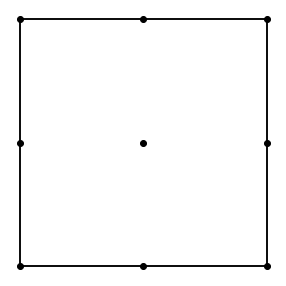}\hspace{1cm}
    \includegraphics[width=.2\textwidth]{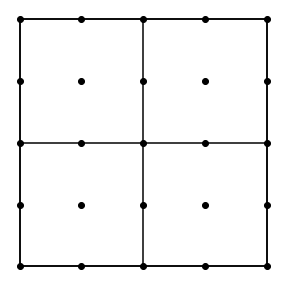}\hspace{1cm}
    \includegraphics[width=.2\textwidth]{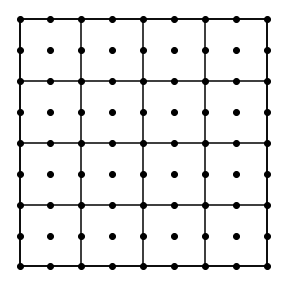}
    \caption{Examples of square domain \( Q^2 \) elements. Rectangles represent the meshes, and dots are the Gauss points within the rectangles. From left to right, the meshes have sizes \( 1 \times 1 \), \( 2 \times 2 \), and \( 4 \times 4 \), each with \( 3 \times 3 \) Gauss points inside. The corresponding values are \( N_x = N_y = 3, 5, 9 \), respectively, and the degrees of freedom are \( N_x \times N_y \).}
    \label{fig:quadrature-grid}
\end{figure}

The finite-dimensional approximation space \( V^h \) is defined as
\begin{equation*}
    V^h = \operatorname{span}\left\{ \psi_i(x)\, \psi_j(y) \;\middle|\; i = 1, \dotsc, N_x;\ j = 1, \dotsc, N_y \right\},
\end{equation*}
where \( \psi_i(x) \) is the \( i \)-th Lagrangian piecewise polynomial of degree \( k \) in the \( x \)-direction, associated with the nodal point \( x_i \). Similarly, \( \psi_j(y) \) is the \( j \)-th Lagrangian piecewise polynomial of degree \( k \)  defined in the \( y \)-direction.

The approximate solution \( u_h(x, y) \in V^h \) can be expressed as
\begin{equation*}
    u_h(x, y) = \sum_{i=1}^{N_x} \sum_{j=1}^{N_y} u_{i,j}\, \psi_i(x)\, \psi_j(y),
\end{equation*}
where \( u_{i,j} = u_h(x_i, y_j) \) due to the interpolation property of the Lagrangian basis functions.

The weak formulation of our PDE \textcolor{black}{with homogeneous boundary condition} is: find \( u_h \in V^h \) such that
\begin{equation}\label{eq:weak-form}
    \int_{\Omega} \nabla u_h \cdot \nabla v_h \, \mathrm{d}x\,\mathrm{d}y + \lambda \int_{\Omega} u_h v_h \, \mathrm{d}x\,\mathrm{d}y + \int_{\Omega} f(x, y, u_h) v_h \, \mathrm{d}x\,\mathrm{d}y = 0, \quad \forall v_h \in V^h.
\end{equation}

The integrals are approximated numerically using Gauss quadrature at the Gauss points \( (x_i, y_j) \). We define the stiffness and mass matrices in the \( x \)- and \( y \)-directions as
\begin{align}
    (S_x)_{i,j} &= \int \frac{\mathrm{d}\psi_i(x)}{\mathrm{d}x} \frac{\mathrm{d}\psi_j(x)}{\mathrm{d}x} \, \mathrm{d}x, &
    (M_x)_{i,j} &= \int \psi_i(x) \psi_j(x) \, \mathrm{d}x, \\
    (S_y)_{i,j} &= \int \frac{\mathrm{d}\psi_i(y)}{\mathrm{d}y} \frac{\mathrm{d}\psi_j(y)}{\mathrm{d}y} \, \mathrm{d}y, &
    (M_y)_{i,j} &= \int \psi_i(y) \psi_j(y) \, \mathrm{d}y,
\end{align}
for \( i, j = 1, \dotsc, N_x \) or \( N_y \), respectively.

The discrete system can then be formulated as
\begin{equation}\label{eq:discrete-system}
    S_x U_h M_y^\top + M_x U_h S_y^\top + \lambda M_x U_h M_y^\top + M_x F_h M_y^\top = 0,
\end{equation}
where \( U_h \) and \( F_h \) are matrices of size \( N_x \times N_y \) defined by
\begin{equation*}
    U_h = \left[ u_{i,j} \right], \quad F_h = \left[ f(x_i, y_j, u_{i,j}) \right].
\end{equation*}

To vectorize the system, we utilize the Kronecker product \( \otimes \) and the property of the vectorization operator \( \operatorname{vec}(\cdot) \):
\begin{equation}\label{eq:vec-property}
    \operatorname{vec}(A X B^\top) = (B \otimes A)\, \operatorname{vec}(X),
\end{equation}
for any given \( A \), \( B \), and \( X \) are matrices of compatible dimensions.

Applying this property to Eq.~\eqref{eq:discrete-system}, we obtain
\begin{equation}\label{eq:vectorized-system}
    \left( M_y \otimes S_x + S_y \otimes M_x + \lambda M_y \otimes M_x \right) \operatorname{vec}(U_h) + \left( M_y \otimes M_x \right) \operatorname{vec}(F_h) = 0.
\end{equation}
Multiplying both sides by \( (M_y \otimes M_x)^{-1} \) the system simplifies to
\begin{equation}\label{eq:final-system}
    \left( I_{N_y} \otimes H_x + H_y \otimes I_{N_x} + \lambda I_{N_y} \otimes I_{N_x} \right) \operatorname{vec}(U_h) + \operatorname{vec}(F_h) = 0,
\end{equation}
where
\begin{equation*}
    H_x = M_x^{-1} S_x, \quad H_y = M_y^{-1} S_y,
\end{equation*}
and \( I_{N_x} \), \( I_{N_y} \) are identity matrices of sizes \( N_x \) and \( N_y \), respectively.

Eq.~\eqref{eq:final-system} represents the final form of the discretized system, suitable for numerical solution.

\textcolor{black}{
For solutions with Dirichlet boundary conditions, i.e., $\alpha_1=0$ in \eqref{eq:main-pde},  consider the approximation space
 \begin{equation*}
    V^h_0 = \operatorname{span}\left\{ v_h\in V^h: v_h|_{\partial \Omega}=0 \right\}.
\end{equation*}
For enforcing the nonhomegeneous Dirichlet boundary condition  $u = \tfrac{\alpha_3}{\alpha_2}$ on $\partial \Omega$,
we use the simple implementation in \cite[Section 6]{li2020superconvergence}, described as follows. The numerical solution
 $u_h \in V^h$ is defined by $u_h=\tilde u_h+u_g$, with $u_g\in V^h$ representing boundary data and $\tilde u_h\in V_0^h$ to be defined later. 
The boundary data $u_g\in V^h$ is the $Q^k$ interpolation of the boundary data function $g(x,y)=\begin{cases}
     0, & (x,y)\notin \partial \Omega \\
    \tfrac{\alpha_3}{\alpha_2}, & (x,y)\in \partial \Omega
 \end{cases}.$ In other words, $u_g\in V_h$ is a piecewise polynomial taking value zero in the interior of $\Omega$, and taking point values of boundary data at the boundary. The part $\tilde u_h\in V_0^h$
 satisfies the weak formulation:  
\begin{equation}\label{eq:weak-form-diri}
    \int_{\Omega} \nabla (\tilde u_h+u_g) \cdot \nabla v_h \, \mathrm{d}x\,\mathrm{d}y + \lambda \int_{\Omega} (\tilde u_h+u_g) v_h \, \mathrm{d}x\,\mathrm{d}y + \int_{\Omega} f(x, y, (\tilde u_h+u_g)) v_h \, \mathrm{d}x\,\mathrm{d}y = 0, \quad \forall v_h\in V_0^h.
\end{equation}
In other words, the numerical solution $u_h\in V^h$ satisfies
\begin{equation*} 
    \int_{\Omega} \nabla u_h  \cdot \nabla v_h \, \mathrm{d}x\,\mathrm{d}y + \lambda \int_{\Omega}  u_h v_h \, \mathrm{d}x\,\mathrm{d}y + \int_{\Omega} f(x, y,  u_h ) v_h \, \mathrm{d}x\,\mathrm{d}y = 0, \quad \forall v_h\in V_0^h.
\end{equation*}   
Such a simple approach allows a finite difference type implementation for treating Dirichlet boundary condition. The Dirichlet boundary data should exist pointwise so that $u_g$ is well defined. 
See \cite[Section 6]{li2020superconvergence} for an {\it a priori} error estimate of this boundary treatment for linear problems.
In the scheme \eqref{eq:weak-form-diri} for Dirichlet boundary conditions, the stiffness matrix can be obtained from taking a submatrix of $S_x$ due to the fact $v_h\in V_0^h$ in  \eqref{eq:weak-form-diri}. See \cite{hu2023positivity,li2020monotonicity,shen2022discrete} for explicit expressions of the matrices for $Q^2$ element and \cite{cross2023monotonicity} for $Q^3$ element.}
\textcolor{black}{For solutions with Neumann boundary conditions, $u_h \in V^h$ satisfies the weak formulation:  
\begin{equation}\label{eq:weak-form-neu}
    \int_{\Omega} \nabla u_h \cdot \nabla v_h \, \mathrm{d}x\,\mathrm{d}y + \lambda \int_{\Omega} u_h v_h \, \mathrm{d}x\,\mathrm{d}y - \int_{\partial \Omega} \frac{\alpha_3}{\alpha_1}  v_h \, \mathrm{d}S+ \int_{\Omega} f(x, y, u_h) v_h \, \mathrm{d}x\,\mathrm{d}y = 0,
\end{equation}
for all $v_h \in V^h$. The integral over $\partial\Omega$ leads to an additional term in the discrete formulation. For instance, if $\Omega = [0,L_x]\times [0,L_y]$ and $\{\psi_i\}$ are the basis functions, then the contribution from the boundary can be written as
\[
  -
  \Bigl[
    \,\psi_{i}(0)\,\int_{0}^{L_y}\!\frac{\alpha_{3}}{\alpha_{1}}\psi_{j}(y)\,dy
    \;+\;
    \psi_{i}(L_x)\,\int_{0}^{L_y}\!\frac{\alpha_{3}}{\alpha_{1}}\psi_{j}(y)\,dy
    \;+\;
    \psi_{j}(0)\,\int_{0}^{L_x}\!\frac{\alpha_{3}}{\alpha_{1}}\psi_{i}(x)\,dx
    \;+\;
    \psi_{j}(L_y)\,\int_{0}^{L_x}\!\frac{\alpha_{3}}{\alpha_{1}}\psi_{i}(x)\,dx
  \Bigr].
\]
}

\begin{remark}\label{Remark_integral}
    \textcolor{black}{Our implementation follows the approach where the $Q^k$ Lagrange basis is defined at Gauss-Lobatto points and integrated using an $k+1$-point Gauss-Lobatto quadrature. In practice, this implementation performs well for cubic nonlinearities in Schrödinger equations, as demonstrated in \cite{chen2024fully}. Furthermore, when using $Q^2$ and $Q^3$ bases, the resulting discrete Laplacian is monotone, which implies certain robustness, e.g., convergence of iterative solvers can be ensured, see \cite[Section 4.3]{zhang-proceeding}. This allows us to get \eqref{eq:discrete-system} even if $F_h$ depends on $u_h$. For example, when $f(u_h)=u_h^3$ with $u_h=\sum_{i}u_i\psi_i$ then we have
    \begin{equation}
        \int_\Omega (u_h)^3 v_h dxdy=\int_\Omega (\sum_{i}u_i\psi_i)^3 v_h dxdy.
    \end{equation} However, under our implementation, we effectively obtain 
        \begin{equation}
        \int_\Omega (u_h)^3 v_h dxdy=\int_\Omega (\sum_{i}u_i^3\psi_i) v_h dxdy
    \end{equation} because or our integration using $k+1$-point Gauss-Lobatto quadrature with $Q^k$ Lagrange basis defined at Gauss-Lobatto points. Moreover, the choice of Lagrange basis at Gauss--Lobatto points ensures that the associated mass matrix is a diagonal mass matrix when we integrate using the same Gauss--Lobatto points because of the Lagrange polynomial property.}
\end{remark}
\subsubsection{Spectral Element Method in 3D}\label{subsec:spectral-element-3d}

Extending the spectral element method to three-dimensional rectangular meshes involves a similar discretization using the \( Q^l \) basis. Let \( (x_i, y_j, z_k) \) denote all the nodal points, where \( i = 1, \dotsc, N_x \); \( j = 1, \dotsc, N_y \); and \( k = 1, \dotsc, N_z \). The finite-dimensional approximation space \( V^h \) is defined as
\begin{equation*}
    V^h = \operatorname{span}\left\{ \psi_i(x)\, \psi_j(y)\, \psi_k(z) \;\middle|\; i = 1, \dotsc, N_x;\ j = 1, \dotsc, N_y;\ k = 1, \dotsc, N_z \right\},
\end{equation*}
where \( \psi_i(x) \), \( \psi_j(y) \), and \( \psi_k(z) \) are the \( i \)-th, \( j \)-th, and \( k \)-th Lagrangian piecewise polynomials of degree \( l \) in the \( x \)-, \( y \)-, and \( z \)-directions, respectively, associated with the nodal points \( x_i \), \( y_j \), and \( z_k \).

The approximate solution \( u_h(x, y, z) \in V^h \) can be expressed as
\begin{equation*}
    u_h(x, y, z) = \sum_{i=1}^{N_x} \sum_{j=1}^{N_y} \sum_{k=1}^{N_z} u_{i,j,k}\, \psi_i(x)\, \psi_j(y)\, \psi_k(z),
\end{equation*}
where \( u_{i,j,k} = u_h(x_i, y_j, z_k) \) due to the interpolation property of the Lagrangian basis functions.

The weak formulation of our PDE remains analogous to the 2D case described in Equation~\eqref{eq:weak-form}, adjusted for three dimensions. Specifically, find \( u_h \in V^h \) such that
\begin{equation}\label{eq:weak-form-3d}
    \int_{\Omega} \nabla u_h \cdot \nabla v_h \, \mathrm{d}x\,\mathrm{d}y\,\mathrm{d}z + \lambda \int_{\Omega} u_h v_h \, \mathrm{d}x\,\mathrm{d}y\,\mathrm{d}z + \int_{\Omega} f(x, y, z, u_h) v_h \, \mathrm{d}x\,\mathrm{d}y\,\mathrm{d}z = 0, \quad \forall v_h \in V^h.
\end{equation}

The integrals are approximated numerically using Gauss quadrature at the Gauss points \( (x_i, y_j, z_k) \). We define the stiffness and mass matrices in the \( x \)-, \( y \)-, and \( z \)-directions as
\begin{align}
    (S_x)_{i,l} &= \int \frac{\mathrm{d}\psi_i(x)}{\mathrm{d}x} \frac{\mathrm{d}\psi_l(x)}{\mathrm{d}x} \, \mathrm{d}x, &
    (M_x)_{i,l} &= \int \psi_i(x) \psi_l(x) \, \mathrm{d}x, \\
    (S_y)_{j,m} &= \int \frac{\mathrm{d}\psi_j(y)}{\mathrm{d}y} \frac{\mathrm{d}\psi_m(y)}{\mathrm{d}y} \, \mathrm{d}y, &
    (M_y)_{j,m} &= \int \psi_j(y) \psi_m(y) \, \mathrm{d}y, \\
    (S_z)_{k,n} &= \int \frac{\mathrm{d}\psi_k(z)}{\mathrm{d}z} \frac{\mathrm{d}\psi_n(z)}{\mathrm{d}z} \, \mathrm{d}z, &
    (M_z)_{k,n} &= \int \psi_k(z) \psi_n(z) \, \mathrm{d}z,
\end{align}
for \( i, l = 1, \dotsc, N_x \); \( j, m = 1, \dotsc, N_y \); and \( k, n = 1, \dotsc, N_z \).

Define \( U_h \) and \( F_h \) as tensors:
\begin{equation*}
    U_h = \left[ u_{i,j,k} \right], \quad F_h = \left[ f(x_i, y_j, z_k, u_{i,j,k}) \right].
\end{equation*}

We obtain the vectorized form from Eq. \eqref{eq:weak-form-3d} by similar computation as 2D:
\begin{equation}\label{eq:vectorized-system-3d}
    \left( M_z \otimes M_y \otimes S_x + M_z \otimes S_y \otimes M_x + S_z \otimes M_y \otimes M_x + \lambda M_z \otimes M_y \otimes M_x \right) \operatorname{vec}(U_h) + \left( M_z \otimes M_y \otimes M_x \right) \operatorname{vec}(F_h) = 0.
\end{equation}

Multiplying both sides by \( (M_z \otimes M_y \otimes M_x)^{-1} \), and defining
\begin{equation*}
    H_x = M_x^{-1} S_x, \quad H_y = M_y^{-1} S_y, \quad H_z = M_z^{-1} S_z,
\end{equation*}
we simplify the system to
\begin{equation}\label{eq:final-system-3d}
    \left( I_{N_z} \otimes I_{N_y} \otimes H_x + I_{N_z} \otimes H_y \otimes I_{N_x} + H_z \otimes I_{N_y} \otimes I_{N_x} + \lambda I_{N_z} \otimes I_{N_y} \otimes I_{N_x} \right) \operatorname{vec}(U_h) + \operatorname{vec}(F_h) = 0,
\end{equation}
where \( I_{N_x} \), \( I_{N_y} \), and \( I_{N_z} \) are identity matrices of sizes \( N_x \), \( N_y \), and \( N_z \), respectively. Eq.~\eqref{eq:final-system-3d} represents the final form of the discretized system in three dimensions. 
\textcolor{black}{We can incorporate the boundary conditions like the 2D case.
}

\subsection{Finite Difference Method}\label{subsec:fdm-2d}

To solve Eq.~\eqref{eq:main-pde}, we employ the finite difference method (FDM) on a uniform grid. The computational domain \( \Omega \) is discretized into \( n_x \) grid points along the \( x \)-axis and \( n_y \) grid points along the \( y \)-axis. The grid spacings are defined as \textcolor{black}{\( h_x = \frac{1}{n_x} \) and \( h_y = \frac{1}{n_y} \)}, with grid points located at
\begin{equation*}
    x_i = i h_x, \quad i = 0, 1, \dotsc n_x,\hbox{~and~}    y_j = j h_y, \quad j = 0, 1, \dotsc, n_y ,\hbox{~when~} \Omega=(0,1)\times(0,1).
\end{equation*}
The numerical solution is denoted by \( (U_h)_{i,j}=u_h(x_i, y_j) \approx u(x_i, y_j) \) at each grid point \( (x_i, y_j) \).

To approximate the second-order derivatives in Eq.~\eqref{eq:main-pde}, we use the central difference scheme. 
%
Then the discretized nonlinear equation at each interior grid point becomes
\begin{equation}\label{eq:fdm-discrete}
    -\left( \frac{(U_h)_{i+1,j} - 2 (U_h)_{i,j} + (U_h)_{i-1,j}}{h_x^2} + \frac{(U_h)_{i,j+1} - 2 (U_h)_{i,j} + (U_h)_{i,j-1}}{h_y^2} \right) + \lambda (U_h)_{i,j} + f(x_i, y_j, (U_h)_{i,j}) = 0.
\end{equation}


We define the following matrix notations:
\begin{itemize}
    \item Solution matrix \( U_h \in \mathbb{R}^{n_x \times n_y} \), where \( (U_h)_{i,j} =u_h(x_i, y_j) \),
    \item Nonlinear matrix \( F_h \in \mathbb{R}^{n_x \times n_y} \), where \( (F_h)_{i,j} = f(x_i, y_j, (U_h)_{i,j}) \),
    \item Mass matrices \( M_x = I_{n_x} \), \( M_y =  I_{n_y} \), where \( I_{n} \) is the identity matrix of size \( n \times n \),
    \item Stiffness matrices \( S_x \in \mathbb{R}^{n_x \times n_x} \) and \( S_y \in \mathbb{R}^{n_y \times n_y} \), defined as tridiagonal matrices with entries
    \begin{equation}
        S_x = \frac{1}{h_x^2}
        \begin{bmatrix}
            2 & -1 & & & \\
            -1 & 2 & -1 & & \\
            & \ddots & \ddots & \ddots & \\
            & & -1 & 2 & -1 \\
            & & & -1 & 2 \\
        \end{bmatrix}
\hbox{~and~}        S_y = \frac{1}{h_y^2}
        \begin{bmatrix}
            2 & -1 & & & \\
            -1 & 2 & -1 & & \\
            & \ddots & \ddots & \ddots & \\
            & & -1 & 2 & -1 \\
            & & & -1 & 2 \\
        \end{bmatrix}. \label{eq:Sy}
    \end{equation}
\end{itemize}

The boundary conditions can be incorporated by modifying the stiffness matrices and the right-hand side appropriately. \textcolor{black}{For problems with Dirichlet boundary conditions, computations are performed exclusively on the interior domain, excluding boundary points (i.e. only using $i\in1,\dots,n_x-1$ and $j\in1,\dots,n_y-1$). Consequently, the dimensionality of the associated matrices is reduced, resulting in a domain of size $\mathbb{R}^{n_x-2 \times n_x-2}, \mathbb{R}^{n_y-2 \times n_y-2}$, respectively. And for $F_h$ part we will change only near the boundary as 
\begin{equation}
(F_h)_{i,j} =
\begin{cases}
    f(x_1, y_j, (U_h)_{1,j}) - \frac{(U_h)_{0,j}}{h_x^2}, & i = 1, \, 2 \leq j \leq n_y - 2, \\
    f(x_{n_x-1}, y_j, (U_h)_{n_x-1,j}) - \frac{(U_h)_{n_x,j}}{h_x^2}, & i = n_x - 1, \, 2 \leq j \leq n_y - 2, \\
    f(x_i, y_1, (U_h)_{i,1}) - \frac{(U_h)_{i,0}}{h_y^2}, & j = 1, \, 2 \leq i \leq n_x - 2, \\
    f(x_i, y_{n_y-1}, (U_h)_{i,n_y-1}) - \frac{(U_h)_{i,n_y}}{h_y^2}, & j = n_y - 1, \, 2 \leq i \leq n_x - 2, \\
    f(x_1, y_1, (U_h)_{1,1}) - \frac{(U_h)_{0,1}}{h_x^2} - \frac{(U_h)_{1,0}}{h_y^2}, & i = 1, \, j = 1, \\
    f(x_1, y_{n_y-1}, (U_h)_{1,n_y-1}) - \frac{(U_h)_{0,n_y-1}}{h_x^2} - \frac{(U_h)_{1,n_y}}{h_y^2}, & i = 1, \, j = n_y - 1, \\
    f(x_{n_x-1}, y_1, (U_h)_{n_x-1,1}) - \frac{(U_h)_{n_x,1}}{h_x^2} - \frac{(U_h)_{n_x-1,0}}{h_y^2}, & i = n_x - 1, \, j = 1, \\
    f(x_{n_x-1}, y_{n_y-1}, (U_h)_{n_x-1,n_y-1}) - \frac{(U_h)_{n_x,n_y-1}}{h_x^2} - \frac{(U_h)_{n_x-1,n_y}}{h_y^2}, & i = n_x - 1, \, j = n_y - 1.
\end{cases}
\end{equation}
For the Neumann boundary condition, $F_h$ change only near the boundary, where ghost points—artificial nodes outside the computational domain—help enforce the condition by approximating derivatives and maintaining numerical consistency. Explicit form is
\begin{equation}
(F_h)_{i,j} =
\begin{cases}
    f(x_0, y_j, (U_h)_{0,j}) - \frac{2\alpha_3}{\alpha_1 h_x}, & i = 0, \, 1 \leq j \leq n_y - 1, \\
    f(x_{n_x}, y_j, (U_h)_{n_x,j}) - \frac{2\alpha_3}{\alpha_1 h_x}, & i = n_x , \, 1 \leq j \leq n_y - 1, \\
    f(x_i, y_0, (U_h)_{i,0}) - \frac{2\alpha_3}{\alpha_1 h_y}, & j = 0, \, 1 \leq i \leq n_x - 1, \\
    f(x_i, y_{n_y}, (U_h)_{i,n_y})- \frac{2\alpha_3}{\alpha_1 h_y}, & j = n_y, \, 1 \leq i \leq n_x - 1, \\
    f(x_0, y_0, (U_h)_{0,0}) - \frac{2\alpha_3}{\alpha_1 h_x} - \frac{2\alpha_3}{\alpha_1 h_y}, & i = 0, \, j = 0, \\
    f(x_0, y_{n_y}, (U_h)_{0,n_y}) - \frac{2\alpha_3}{\alpha_1 h_x} - \frac{2\alpha_3}{\alpha_1 h_y}, & i = 0, \, j = n_y, \\
    f(x_{n_x}, y_0, (U_h)_{n_x,0}) -  \frac{2\alpha_3}{\alpha_1 h_x} - \frac{2\alpha_3}{\alpha_1 h_y}, & i = n_x, \, j = 0, \\
    f(x_{n_x}, y_{n_y}, (U_h)_{n_x,n_y}) - \frac{2\alpha_3}{\alpha_1 h_x} - \frac{2\alpha_3}{\alpha_1 h_y}, & i = n_x, \, j = n_y,
\end{cases}
\end{equation}
with a different stiffness matrices \( S_x \in \mathbb{R}^{n_x \times n_x} \) and \( S_y \in \mathbb{R}^{n_y \times n_y} \), defined as 
\begin{equation}
        S_x = \frac{1}{h_x^2}
        \begin{bmatrix}
            2 & -2 & & & \\
            -1 & 2 & -1 & & \\
            & \ddots & \ddots & \ddots & \\
            & & -1 & 2 & -1 \\
            & & & -2 & 2 \\
        \end{bmatrix}
\hbox{~and~}        S_y = \frac{1}{h_y^2}
        \begin{bmatrix}
            2 & -2 & & & \\
            -1 & 2 & -1 & & \\
            & \ddots & \ddots & \ddots & \\
            & & -1 & 2 & -1 \\
            & & & -2 & 2 \\
        \end{bmatrix}.
    \end{equation}}

    The discretized system can be compactly written as
\begin{equation}\label{eq:fdm-matrix}
    S_x U_h M_y^\top + M_x U_h S_y^\top + \lambda M_x U_h M_y^\top + M_x F_h M_y^\top = 0.
\end{equation}

Alternatively, the system can be vectorized using the Kronecker product \( \otimes \) and the \(\operatorname{vec}(\cdot)\) operator:
\begin{equation}\label{eq:fdm-vectorized_2}
    \left( M_y \otimes S_x + S_y \otimes M_x + \lambda M_y \otimes M_x \right) \operatorname{vec}(U_h) + \left( M_y \otimes M_x \right) \operatorname{vec}(F_h) = 0.
\end{equation}

Multiplying both sides by \( (M_y \otimes M_x)^{-1} \) and defining \( H_x = M_x^{-1} S_x \), \( H_y = M_y^{-1} S_y \), we obtain
\begin{equation}\label{eq:fdm-final}
    \left( I_{n_y} \otimes H_x + H_y \otimes I_{n_x} + \lambda I_{n_y} \otimes I_{n_x} \right) \operatorname{vec}(U_h) + \operatorname{vec}(F_h) = 0,
\end{equation}
where \( I_{n_x} \) and \( I_{n_y} \) are identity matrices of sizes \( n_x \times n_x \) and \( n_y \times n_y \), respectively.

Equation~\eqref{eq:fdm-final} represents the final form of the discretized nonlinear system, suitable for numerical solutions using iterative methods designed for large sparse systems. We can also do similar computations for 3D cases and the result will be similar to Eq. \eqref{eq:final-system-3d}.

For the 3D case, we can get the system as
\begin{equation}\label{eq:final-system-3d-fdm}
    \left( I_{n_z} \otimes I_{n_y} \otimes H_x + I_{n_z} \otimes H_y \otimes I_{n_x} + H_z \otimes I_{n_y} \otimes I_{n_x} + \lambda I_{n_z} \otimes I_{N_y} \otimes I_{n_x} \right) \operatorname{vec}(U_h) + \operatorname{vec}(F_h) = 0,
\end{equation}
where \( I_{n_x} \), \( I_{n_y} \), and \( I_{n_z} \) are identity matrices of sizes \( n_x \), \( n_y \), and \( n_z \), respectively with
\begin{equation*}
    H_x = S_x, \quad H_y =  S_y, \quad H_z = S_z.
\end{equation*}
Remark that we know the $S_x, S_y$ and $S_z$ as
 \begin{equation*}
        S_x = \frac{1}{h_x^2}
        \begin{bmatrix}
             2 & -1 & \\
            \ddots & \ddots & \ddots \\
            &  -1 & 2\\
        \end{bmatrix}
\hbox{~,~}        S_y = \frac{1}{h_y^2}
        \begin{bmatrix}
             2 & -1 & \\
            \ddots & \ddots & \ddots \\
            &  -1 & 2\\
    \end{bmatrix}\hbox{~and~}        S_z = \frac{1}{h_z^2}
        \begin{bmatrix}
             2 & -1 & \\
            \ddots & \ddots & \ddots \\
            &  -1 & 2\\
        \end{bmatrix}.
    \end{equation*}

\textcolor{black}{For the boundary conditions, we can perform similar computations as in the 2D case.
}

\section{The Quasi-Newton Method}\label{IM}
In this section, we introduce the Quasi-Newton method. We modify the nonlinear term by replacing it with \( \beta_n I \), where \( I \) is the identity matrix and \( \beta_n \) is a constant. This modification allows us to fully exploit the structure of tensor products and Kronecker products. The primary reason for the increased computational efficiency of this approach is discussed in \S~\ref{implementation}.

We rewrite the discretized nonlinear system as follows:
\begin{equation}\label{eq:F_u}
    \textcolor{black}{F^h(U_h) = A_h\operatorname{vec}(U_h) + N_h(\vec{x}, U_h),}
\end{equation}
where \( A_h \) is a matrix that describes linear operator, \( N_h(\vec{x}, U_h) \) represents the nonlinear component, and \( U_h \) is the numerical solution. From now on
\(\textcolor{black}{\operatorname{vec}(U_h)} \) is the vectorization of the solution matrix \( (U_h)_{i,j} \) or tensor \( (U_h)_{i,j,k} \), and \( N_h(\vec{x}, U_h) = \operatorname{vec}([f(x_i,y_j,(U_h)_{i,j})]) \) in 2D or \( N_h(\vec{x}, U_h) = \operatorname{vec}([f(x_i,y_j,z_k,(U_h)_{i,j,k})]) \) in 3D is the vectorization of the nonlinear function evaluated at each grid point as in \S \ref{Discret_2d}. The matrix \( A_h \) is given by
\begin{equation}\label{eq:A_definition_2D}
    A_h = I_{y} \otimes H_x + H_y \otimes I_{x} + \lambda I_{y} \otimes I_{x} \hbox{~in the 2D case}
\end{equation}
and 
\begin{equation}\label{eq:A_definition_3D}
    A_h = I_{z} \otimes I_{y} \otimes H_x + I_{z} \otimes H_y \otimes I_{x} + H_z \otimes I_{y} \otimes I_{x} + \lambda I_{z} \otimes I_{y} \otimes I_{x} \hbox{~in the 3D case}
\end{equation}

Newton's method  for solving Eq.~\eqref{eq:F_u} reads as 
\begin{equation}\label{eq:newton}
    \textcolor{black}{\operatorname{vec}(U_{n+1}) = \operatorname{vec}(U_{n}) - [DF^h(U_n)]^{-1} F^h(U_n),}
\end{equation}
where \( DF^h(\cdot)=A_h+(N_h)_U(\cdot) \) is the Jacobian matrix of \( F^h \). 
See \cite[Theorem 8]{zhang-proceeding} for convergence of Newton's method using $Q^2$ spectral element method and second order finite difference for a quasi-linear PDE.

The computational bottleneck in Newton's method is the inversion of the Jacobian matrix at each iteration. To take advantage of tensor product structure,  we approximate the Jacobian matrix  and define the quai-Newton method as:
\begin{equation}\label{eq:identity-method}
    \textcolor{black}{\operatorname{vec}(U_{n+1}) = \operatorname{vec}(U_{n}) - [A_h + \beta_n I]^{-1} F^h(U_n),}
\end{equation}
where \( I \) is the identity matrix, and \( A_h + \beta_n I \) serves as an approximation to the Jacobian \( DF^h(u_n) \) by choosing \begin{equation*}
    \color{black}\beta_n = \frac{\lambda_{\max}(\textcolor{black}{(N_h)_U(U_n)}) + \lambda_{\min}(\textcolor{black}{(N_h)_U(U_n)})}{2}.
\end{equation*}

\subsection{Implementation with tensor product}\label{implementation}



We illustrate the procedure of computing the inverse $A_h+\beta_n I$ with the tensor product, namely,
\begin{equation}
    \left( I_y \otimes H_x + H_y \otimes I_x + (\lambda+\beta_n) I_y \otimes I_x \right) \operatorname{vec}(U) = \operatorname{vec}(F^h(U_n)),\label{eq:linear_system}
\end{equation}
By utilizing the diagonalizations of \( H_x \) and \( H_y \) as in \cite{liu2024simple} \textcolor{black}{and} \S \ref{Diagonal}, we have
\begin{equation*}
    H_x = T_x \Lambda_x T_x^{-1}, \quad H_y = T_y \Lambda_y T_y^{-1}.\nonumber
\end{equation*}
Thus Equation~\eqref{eq:linear_system} becomes
\begin{equation*}
    \left( T_y \otimes T_x \right) \left( I_y \otimes \Lambda_x + \Lambda_y \otimes I_x + (\lambda+\beta_n) I_y \otimes I_x \right) \left( T_y^{-1} \otimes T_x^{-1} \right) \operatorname{vec}(U) = \operatorname{vec}(F^h(U_n)).
\end{equation*}

Multiplying both sides by \( \left( T_y^{-1} \otimes T_x^{-1} \right) \), we obtain
\begin{equation*}
    \left( I_y \otimes \Lambda_x + \Lambda_y \otimes I_x + (\lambda+\beta_n) I_y \otimes I_x \right) \operatorname{vec}\left( T_x^{-1} U T_y^{-\top} \right) = \operatorname{vec}\left( T_x^{-1} F^h(U_n) T_y^{-\top} \right).
\end{equation*}
Since \( I_y \otimes \Lambda_x + \Lambda_y \otimes I_x + \lambda I_y \otimes I_x \) is diagonal, its inverse is straightforward to compute element-wise. Therefore, the solution \( U \) can be expressed as
\begin{equation}\label{fin_inverse}
    U = T_x \left( \left( T_x^{-1} F^h(U_n) T_y^{-\top} \right) \oslash \left( (\Lambda_x)_{i,i} + (\Lambda_y)_{j,j} + \lambda+\beta_n \right) \right) T_y^\top,
\end{equation}
where \( \oslash \) denotes element-wise division. 
This technique can be extended to three dimensions, which significantly improves computational speed as the dimensionality of the domain increases \cite{liu2024simple}.

The primary reason for the reduction in computational time is that we do not need to compute the matrix inverse at every iteration, we only need to compute the inverse once. The commutativity between $A_h$ and $I$ allows us to diagonalize the matrices, enabling us to compute the inverse efficiently. Specifically, once the matrices are diagonalized, the inverse computation reduces to inverting diagonal matrices, which is straightforward.

Moreover, this approach is faster than performing full matrix multiplication because we are operating on significantly smaller matrices. After diagonalization, we only need to perform element-wise division, which is computationally inexpensive compared to matrix inversion or multiplication.

By leveraging the diagonalization and commutativity of the matrices, we significantly reduce the computational cost, thereby improving the overall efficiency of the method.

\textcolor{black}{To efficiently compute the tensor transformation, we leverage GPU-based operations to accelerate large-scale computations. The MATLAB implementation for the transformation involving the Kronecker product and vectorization is given by:
\begin{equation}
    \operatorname{vec}(A X B^\top) = (B \otimes A)\, \operatorname{vec}(X),
\end{equation}
where \( B \) and \( A \) are tensors, and \( X \) is a tensor whose vectorization is denoted by \( \operatorname{vec}(X) \).
In MATLAB, this operation can be implemented using the following code:}

\begin{center}
\begin{minipage}{0.6\textwidth}\color{black}
\begin{verbatim}
% Given Tensor A, B, and X
    X = gpuArray(X);  % Move data to GPU
    A = gpuArray(A);
    B = gpuArray(B);
% Equivalent to matrix muliplication X(B transpose).
u1 = tensorprod(X, B', 2, 1);

% Equivalent to matrix muliplication AX(B transpose)
u1 = squeeze(tensorprod(A, u1, 2, 1));
\end{verbatim}
\end{minipage}
\end{center}
\textcolor{black}{The key advantage of this approach is the use of parallel computation on the GPU, which is particularly beneficial for large tensors.}

\textcolor{black}{For the 3D case, the MATLAB implementation can be written as follows:
\begin{equation}
    (C \otimes B \otimes A)\, \operatorname{vec}(X),
\end{equation}
this transformation can be implemented as:}
\begin{center}
\begin{minipage}{0.6\textwidth}\color{black}
\begin{verbatim}
% Given Tensor A, B, C, and X
% Given Tensor A, B, and X
    X = gpuArray(X);  % Move data to GPU
    A = gpuArray(A);
    B = gpuArray(B);
    C = gpuArray(C);
u1 = tensorprod(X, C', 3, 1);
u1 = pagemtimes(u1, B');
u1 = squeeze(tensorprod(A, u1, 2, 1));
\end{verbatim}
\end{minipage}
\end{center}

\subsection{Diagonalization of Matrices}\label{Diagonal}

We observe that the mass matrix \( M_x \) is diagonal with positive weights, hence it is positive definite. We can express the matrix \( H_x \) as:
\begin{equation*}
    H_x = M_x^{-1} S_x = M_x^{-1/2} \left( M_x^{-1/2} S_x M_x^{-1/2} \right) M_x^{1/2}.
\end{equation*}
Since \( M_x^{-1/2} S_x M_x^{-1/2} \) is symmetric, it can be diagonalized:
\begin{equation*}
    M_x^{-1/2} S_x M_x^{-1/2} = \tilde{T}_x \Lambda_x \tilde{T}_x^{-1},
\end{equation*}
where \( \tilde{T}_x \) is the matrix of eigenvectors and \( \Lambda_x \) is the diagonal matrix of eigenvalues.

Defining
\begin{equation*}
    T_x = M_x^{-1/2} \tilde{T}_x,
\end{equation*}
we can write
\begin{equation*}
    H_x = M_x^{-1} S_x = T_x \Lambda_x T_x^{-1}.
\end{equation*}
Thus, \( H_x \) is diagonalizable.

\subsection{Extend to nonlinear
systems of differential equations}\label{ext_nonli_sys}

When dealing with systems of differential equations, we may need to compute the eigenvalues of non-diagonal matrices. In this sub-section, we consider the following system of two differential equations which can be applied to $n$ equations straightforwardly:
\begin{equation}\label{eq:system}
\begin{cases}
-\Delta u + \lambda_1 u + f(x, u, v) = 0, \quad \text{in } \Omega \hbox{~and~}
\alpha_1 \frac{\partial u}{\partial n}+\beta_1 u=\gamma_1, \quad \text{on } \partial\Omega,
\\
  -\Delta v + \lambda_2 v + g(x, u, v) = 0, \quad \text{in } \Omega \hbox{~and~}
\alpha_2 \frac{\partial v}{\partial n}+\beta_2 v=\gamma_2, \quad \text{on } \partial\Omega.
\end{cases}
\end{equation}
Then the discretized system is written as 
\begin{equation}
    F^h(U,V)=\begin{bmatrix}
        A_1^h & 0 \\
        0 & A_2^h
    \end{bmatrix}(\textcolor{black}{\operatorname{vec}(U),\operatorname{vec}(V)})^\top +  N_h(\vec{x},U,V)
\end{equation}
where
\begin{equation*}
    \begin{cases}
        A_1^h=I_{y} \otimes H_x + H_y \otimes I_{x} + \lambda_1 I_{y} \otimes I_{x},
        \\
        A_2^h=I_{y} \otimes H_x + H_y \otimes I_{x} + \lambda_2 I_{y} \otimes I_{x},
    \end{cases}
    \hbox{~in the 2D case}
\end{equation*}
and
\begin{equation*}
    \begin{cases}
        N_1=vec([f(x_i,y_j,U_{i,j},V_{i,j})]),
        \\
        N_2=vec([g(x_i,y_j,U_{i,j},V_{i,j})]),
    \end{cases} \hbox{~in the 2D case.}
\end{equation*}
Next, we compute the eigenvalue of  the following block matrix:
\begin{equation*}
   DN_h=\begin{bmatrix}
        (N_1)_U & (N_1)_V \\
        (N_2)_U & (N_2)_V
    \end{bmatrix}.
\end{equation*}
where, in our specific case, \( (N_1)_U, (N_1)_V, (N_2)_U, \) and \( (N_2)_V \) are diagonal matrices. 




To compute the eigenvalues \( \lambda \) of \( DN_h \), we have
\begin{equation*}
    \det\left( DN_h - \lambda I \right) = \det\left( ((N_1)_U - \lambda I)((N_2)_V - \lambda I) - (N_1)_V (N_2)_U \right).
\end{equation*}
Since  \( (N_1)_U, (N_1)_V, (N_2)_U, \) and \( (N_2)_V \) are diagonal matrices, their products and differences remain diagonal, and thus the determinant simplifies to:
\begin{equation*}
    \det\left( DN_h - \lambda I \right) = \prod_{i=1}^n \left(  (((N_1)_U)_{i,i} - \lambda)( ((N_2)_V)_{i,i} - \lambda) - ((N_1)_V)_{i,i} ((N_2)_U)_{i,i} \right).
\end{equation*}
To compute the eigenvalues \( \lambda \), we solve the quadratic equations:
\begin{equation}\label{eq:eigen_quad}
    (((N_1)_U)_{i,i} - \lambda)( ((N_2)_V)_{i,i} - \lambda) - ((N_1)_V)_{i,i} ((N_2)_U)_{i,i} = 0, \quad \text{for } i = 1, 2, \dots, n.
\end{equation}
Thus we have:
\begin{equation}\label{eq:eigenvalues}
    \lambda_i = \frac{((N_1)_U)_{i,i} + ((N_2)_V)_{i,i} \pm \sqrt{(((N_1)_U)_{i,i} - ((N_2)_V)_{i,i})^2 + 4 ((N_1)_V)_{i,i} ((N_2)_U)_{i,i}}}{2}, \quad \text{for } i = 1, 2, \dots, n.
\end{equation}
Thus, even when dealing with a system of equations where the nonlinear term leads to a non-diagonal \( N \), we can compute the eigenvalues by using the block matrix \( DN_h \) and solving for \( \lambda \) using the determinant of \( DN_h - \lambda I \). This allows us to proceed with the optimization of \( \beta \) in our method, ensuring convergence and efficiency.

\subsection{Eigenvalue Problem}\label{explanation_eigp}

We consider the eigenvalue problem:
\begin{equation}\label{Eigenvalue_p}
\begin{cases} 
    -\Delta u + f(u) = \lambda u, & \text{in } \Omega, \\ 
    \partial_n u = 0, & \text{on } \partial \Omega,
\end{cases}
\end{equation}
where both the eigenfunction \( u \) and the eigenvalue \( \lambda \) are unknown. Direct application of Newton's method is not straightforward in this setting due to the presence of the unknown parameter \( \lambda \). To address this issue, we employ an iterative approach by expressing \( \lambda \) in terms of \( u \):
\begin{equation*}
    \lambda = \frac{\int_{\Omega}(-\Delta u + f(u))u}{\int_\Omega (u)(u)},
\end{equation*}
and substituting this expression back into equation~\eqref{Eigenvalue_p}. This reformulation allows us to apply Newton's method to solve for \( u \) while updating \( \lambda \) accordingly. Furthermore, in \S~\ref{3d_eigen_p}, we include a normalization step for \( u \); after each application of Quasi-Newton's method, we normalize \( u \) to control its magnitude. More precisely, we define the discrete nonlinear operator \( F^h(U_h) \) as
\begin{equation*}
    F^h(U_h) = A_h \textcolor{black}{\operatorname{vec}(U_h)} + N_h(U_h),
\end{equation*}
where \( \textcolor{black}{\operatorname{vec}(U_h)} \) is the discrete solution vector, and the matrices \( A_h \) and \( N_h(U_h) \) are defined as follows.

In the \textcolor{black}{three-}dimensional case, the matrix \( A_h \) is given by
\begin{equation*}
    \color{black}A_h = I_{z} \otimes I_{y} \otimes H_x + I_{z} \otimes H_y \otimes I_{x} + H_z \otimes I_{y} \otimes I_{x},
\end{equation*}
where \( H_x \) and \( H_y \) are discretization matrices in the \( x \)- and \( y \)-directions, respectively. The nonlinear term \( N_h(U_n) \) is defined as
\begin{equation*}
\color{black}\lambda_n=\frac{\left<A_h{\operatorname{vec}(U_n)}+f({\operatorname{vec}(U_n)}),{\operatorname{vec}(U_n)}\right>}{\left<{\operatorname{vec}(U_n)},{\operatorname{vec}(U_n)}\right>}
\end{equation*}
\begin{equation*}
    \color{black}N_h(U_n)=f({\operatorname{vec}(U_n)})-\lambda_n{\operatorname{vec}(U_n)},
\end{equation*}
with \(\left<\textcolor{black}{\operatorname{vec}(U^1_h)},\textcolor{black}{\operatorname{vec}(U^2_h)}\right>=\operatorname{vec}(U^1_h)^\top \left( M_z \otimes M_y \otimes M_x \right) \operatorname{vec}(U^2_h)\) for 3D. The iterative scheme is given by
\begin{equation}\label{system_qusi}
    \textcolor{black}{\operatorname{vec}(\Tilde{U}_{n+1})} = \textcolor{black}{\operatorname{vec}({U}_{n})} - [A_h + \beta_n I]^{-1} F^h(U_n),
\end{equation}
followed by the normalization\textcolor{black}{\begin{equation}\label{system_normal}
    \operatorname{vec}({U_{n+1}}) =\frac{\operatorname{vec}(\Tilde{U}_{n+1})}{\|\operatorname{vec}(\Tilde{U}_{n+1})\|_{L^2}}= \frac{\operatorname{vec}(\Tilde{U}_{n+1})}{\sqrt{\left<\operatorname{vec}(\Tilde{U}_{n+1}),\operatorname{vec}(\Tilde{U}_{n+1})\right>}}.
\end{equation}}
We \textcolor{black}{start from $U_0$ with $\left<\textcolor{black}{\operatorname{vec}(U_0)},\textcolor{black}{\operatorname{vec}(U_0)}\right>=1$ and} repeat the iterative steps defined by Eq. \eqref{system_qusi}, \eqref{system_normal} until the residual \( F^h(U_n) \) is sufficiently small. The results obtained using this method are presented in \S~\ref{3d_eigen_p}.

\section{Convergence Analysis}\label{Convergence}
\textcolor{black}{Here we use matrix norm as
\begin{equation*}
\|A\|=\sup_{\|U\|_{L^2}=1}\frac{\|A \operatorname{vec}(U)\|_{L^2}}{\|U\|_{L^2}}, \text{~with~} \|U\|_{L^2}=
\begin{cases}
    \operatorname{vec}(U)^\top \left( M_z \otimes M_y \otimes M_x \right) \operatorname{vec}(U) \text{~when~}U \text{~is 3 dimension tensor}\\
    \operatorname{vec}(U)^\top \left( M_y \otimes M_x \right) \operatorname{vec}(U) \text{~when~}U \text{~is 2 dimension tensor}\\
\end{cases}
\end{equation*}
}
\begin{thm}\label{thm:convergence}
Let \( F^h(U) = A_h \textcolor{black}{\operatorname{vec}(U)} + N_h(\vec{x}, U) \), where \textcolor{black}{\( A_h \)} is a symmetric positive definite (SPD) matrix, and the Jacobian of \textcolor{black}{\( N_h = \operatorname{vec}([f(\vec{x},U)])  \) with respect to \( U \) denotes \( (N_h)_U(\vec{x}, U)  \) where $f(\vec{x},~\cdot~)\in C^2(\mathbb{R})$}. Suppose that \( (N_h)_U \) is SPD near the solution \( U^\ast \) (\( U^\ast \) be the exact solution satisfying \( F^h(U^\ast) = 0 \)). Then, the iterative scheme
\begin{equation}\label{eq:identity_method}
    \textcolor{black}{\operatorname{vec}(U_{n+1}) = \operatorname{vec}(U_n) - [A_h + \beta_n I]^{-1} F^h(U_n),}
\end{equation}
will converge locally to \( U^\ast \) \textcolor{black}{with $\|U^\ast-U_0\|_{L^2}\ll1$} when
\begin{equation}
    \color{black}\beta_n = \frac{\lambda_{\max}(\textcolor{black}{(N_h)_U(U_n)}) + \lambda_{\min}(\textcolor{black}{(N_h)_U(U_n)})}{2},
\end{equation}
where \textcolor{black}{\( \lambda_{\max}((N_h)_U(U_n)) \) and \( \lambda_{\min}((N_h)_U(U_n)) \)} are the maximum and minimum eigenvalues of \textcolor{black}{\( (N_h)_U (U_n)\)}.
\end{thm}

\begin{proof}
By defining the error at iteration \( n \) as $   e_n = \operatorname{vec}(U^\ast) - \operatorname{vec}(U_n)$,
we obtain:
\begin{equation}\label{eq:error-recursion}
    e_{n+1} = e_n + [A_h + \beta_n I]^{-1} F^h(U_n).
\end{equation}

Using a Taylor expansion of \( F^h(U) \) at \( U^\ast \), \textcolor{black}{because $f(x,~\cdot~)\in C^2(\mathbb{R})$} we get:
\begin{equation}\label{eq:F_u_n}
    \color{black} F^h(U^\ast)= 0 = F^h(U_n) + DF^h(U_n) e_n + R_n,
\end{equation}
where \( R_n = O(\|e_n\|^2_{L^2}) \) represents the higher-order remainder terms \textcolor{black}{which will depends on $U^*$ locally}. 

Substituting Equation~\eqref{eq:F_u_n} into Equation~\eqref{eq:error-recursion}, we have:
\begin{align}
    e_{n+1} &= e_n - [A_h + \beta_n I]^{-1} \left( DF^h\textcolor{black}{(U_n)} e_n + R_n \right) \nonumber \\
    &= \left( I - [A_h + \beta_n I]^{-1} DF^h\textcolor{black}{(U_n)} \right) e_n - [A_h + \beta_n I]^{-1} R_n. \label{eq:error-update}
\end{align}

Taking norms on both sides of Equation~\eqref{eq:error-update}, we obtain:
\begin{equation}\label{eq:error-norm}
    \|e_{n+1}\|_{L^2} \leq \left\| I - [A_h + \beta_n I]^{-1} DF^h\textcolor{black}{(U_n)} \right\| \|e_n\|_{L^2} + \left\| [A_h + \beta_n I]^{-1} \right\| \|R_n\|_{L^2}.
\end{equation}

Since \( R_n = O(\|e_n\|^2_{L^2}) \) represents the higher-order remainder terms, there exists a constant \( \tau > 0 \) such that \( \|R_n\|_{L^2} \leq \tau \|e_n\|^2_{L^2} \) when \( \|e_n\|_{L^2} \) is sufficiently small. Therefore, Equation~\eqref{eq:error-norm} becomes:
\begin{equation}\label{eq:error-norm-final}
    \|e_{n+1}\|_{L^2} \leq \left( \left\| I - [A_h + \beta_n I]^{-1} DF^h\textcolor{black}{(U_n)} \right\| + \tau \left\| [A_h + \beta_n I]^{-1} \right\| \|e_n\|_{L^2} \right) \|e_n\|_{L^2}.
\end{equation}

Next, we will prove
\begin{equation}\label{eq:convergence-condition}
    \rho \left( I - [A_h + \beta_n I]^{-1} DF^h\textcolor{black}{(U_n)} \right) < 1,
\end{equation}
where \( \rho(\cdot) \) denotes the spectral radius of a matrix.

The Jacobian matrix \( DF^h\textcolor{black}{(U_n)} \) is the sum of two matrices:

By denoting $
DF^h(U_n) = A_h + \textcolor{black}{(N_h)_U(U_n)},
$
%
%
we have
\begin{align*}
I - [A_h + \beta_n I]^{-1} (A_h + \textcolor{black}{(N_h)_U(U_n)}) &= [A_h + \beta_n I]^{-1} \left( \beta_n I - \textcolor{black}{(N_h)_U(U_n)}\right).
\end{align*}
%
Moreover, we have
\begin{equation*}
\left\| [A_h + \beta_n I]^{-1} (\beta_n I - \textcolor{black}{(N_h)_U(U_n)}) \right\| \leq \left\| [A_h + \beta_n I]^{-1} \right\| \left\| \beta_n I - \textcolor{black}{(N_h)_U(U_n)} \right\|.
\end{equation*}

Let \( \Lambda(A_h) \) and \( \Lambda(\textcolor{black}{(N_h)_U(U_n)}) \) denote the spectra (sets of eigenvalues) of \( A_h \) and \( \textcolor{black}{(N_h)_U(U_n)} \), respectively.

Since \( A_h + \beta_n I \) is invertible for \( \beta_n > -\lambda_{\min}(A_h) \), where \( \lambda_{\min}(A_h) \) is the smallest eigenvalue of \( A_h \), we have:

\[
\left\| [A_h + \beta_n I]^{-1} \right\| = \frac{1}{\min_{\lambda \in \Lambda(A_h)} |\lambda + \beta_n|}.
\]

Since
\[
\left\| \beta_n I - \textcolor{black}{(N_h)_U(U_n)} \right\| = \max_{\mu \in \Lambda(\textcolor{black}{(N_h)_U(U_n)})} |\beta_n - \mu|,
\]
we obtain:
\begin{equation}\label{eq:norm_ratio}
\left\| [A_h + \beta_n I]^{-1} (\beta_n I - \textcolor{black}{(N_h)_U(U_n)}) \right\| \leq \frac{\max_{\mu \in \Lambda(\textcolor{black}{(N_h)_U(U_n)})} |\beta_n - \mu|}{\min_{\lambda \in \Lambda(A_h)} |\lambda + \beta_n|}.
\end{equation}

By choosing
\[
\textcolor{black}{\beta^\ast := \frac{\lambda_{\min}((N_h)_U(U_n))+ \lambda_{\max}((N_h)_U(U_n))}{2},}
\]
we have:

\[
\max_{\mu \in \Lambda(\textcolor{black}{(N_h)_U(U_n)})} |\beta^\ast - \mu| = -\lambda_{\min}(\textcolor{black}{(N_h)_U(U_n)})+\frac{\lambda_{\max}(\textcolor{black}{(N_h)_U(U_n)}) + \lambda_{\min}(\textcolor{black}{(N_h)_U(U_n)})}{2},
\]

\[
\min_{\lambda \in \Lambda(A_h)} |\lambda + \beta^\ast| = \lambda_{\min}(A_h) + \beta^\ast = \lambda_{\min}(A_h)+ \frac{\lambda_{\min}(\textcolor{black}{(N_h)_U(U_n)}) + \lambda_{\max}(\textcolor{black}{(N_h)_U(U_n)})}{2}.
\]

Since \(\beta^\ast>0 \) (because \( \textcolor{black}{(N_h)_U(U_n)} \) is SPD near \( U^\ast \)), we  ensure that:
\begin{equation}\label{inequality}
\frac{\max_{\mu \in \Lambda(\textcolor{black}{(N_h)_U(U_n)})} |\beta^\ast - \mu|}{\min_{\lambda \in \Lambda(A_h)} |\lambda + \beta^\ast|} = \frac{-\lambda_{\min}(\textcolor{black}{(N_h)_U(U_n)})+\frac{\lambda_{\max}(\textcolor{black}{(N_h)_U(U_n)}) + \lambda_{\min}(\textcolor{black}{(N_h)_U(U_n)})}{2}}{\lambda_{\min}(A_h) + \frac{\lambda_{\min}(\textcolor{black}{(N_h)_U(U_n)}) +\lambda_{\max}(\textcolor{black}{(N_h)_U(U_n)})}{2}} < 1.
\end{equation}

Since \( \lambda_{\min}(A_h) > 0 \) (because \( A_h \) is SPD), \( \lambda_{\min}(\textcolor{black}{(N_h)_U(U_n)}) > -\lambda_{\min}(A_h) \), the inequality is satisfied. Therefore, choosing \( \beta_n = \frac{\lambda_{\min}((N_h)_U(U)_n)) + \lambda_{\max}((N_h)_U(U)_n)}{2} \) ensures that \( \rho\left( I - [A_h + \textcolor{black}{\beta_n} I]^{-1} (A_h + (N_h)_U) \right) < 1 \). \textcolor{black}{Assuming that $\|U^\ast - U_0\| \leq \epsilon \ll 1$, from Equations \eqref{eq:error-norm-final} and \eqref{eq:norm_ratio}, we have
\begin{equation}
    \|e_1\|_{L^2} \leq \frac{\max_{\mu \in \Lambda((N_h)_U)} |\beta_0 - \mu| + \gamma \epsilon}{\min_{\lambda \in \Lambda(A_h)} |\lambda + \beta_0|} \|e_0\|_{L^2} < \|e_0\|_{L^2} \leq \epsilon.
\end{equation}
Inductively, we get
\begin{equation}
    \|e_n\|_{L^2} \leq \prod_{k=0}^{n-1}\left( \frac{\max_{\mu \in \Lambda((N_h)_U(U_k))} |\beta_k - \mu| + \gamma \epsilon}{\min_{\lambda \in \Lambda(A_h)} |\lambda + \beta_k|} \right) \|e_0\|_{L^2} \leq \epsilon.
\end{equation}
Using a Taylor expansion, we have
\begin{equation}
    F^h(U_n) = 0 = F^h(U^\ast) - DF^h(U^\ast) e_n + R_n,
\end{equation}
and from Equation \eqref{eq:F_u_n}, we obtain
\begin{equation}
    \|DF^h(U^\ast) - DF^h(U_n)\| = O(\|e_n\|).
\end{equation}
This leads to
\begin{equation}
    \frac{\lambda_{\min}((N_h)_U(U_n)) + \lambda_{\max}((N_h)_U(U_n))}{2} < C < \infty,
\end{equation}
for some constant $C$, with Equation \eqref{inequality} implying
\begin{equation}
    \|e_n\|_{L^2} \leq \prod_{k=0}^{n-1}\left( \frac{C + \gamma \epsilon}{\lambda_{\min}(A_h) + C} \right) \|e_0\|_{L^2} \to_{n \to \infty} 0,
\end{equation}
 leading to convergence of the quasi-Newton method. 
}

\end{proof}

\begin{remark} {\bf Optimal $\beta$:}
    \textcolor{black}{Assuming \( A_h \) is symmetric positive definite (SPD) and \( -\lambda_{\min}(A_h) < \lambda_{\min}((N_h)_U(U_n)) \). We can define
\[
g(\beta) := \frac{\max_{\mu \in \Lambda((N_h)_U(U_n))} |\beta - \mu|}{\min_{\lambda \in \Lambda(A_h)} |\lambda + \beta|}
\]
    and this can expressed with piecewise function:}

\[
g(\beta) =
\begin{cases}
\displaystyle \frac{\lambda_{\max}(\textcolor{black}{(N_h)_U(U_n)}) - \beta}{\beta + \lambda_{\min}(A_h)}, & \text{if } -\lambda_{\min}(A_h)< \beta < \frac{\lambda_{\min}(\textcolor{black}{(N_h)_U(U_n)}) + \lambda_{\max}(\textcolor{black}{(N_h)_U(U_n)})}{2}, \\
\displaystyle \frac{\beta - \lambda_{\min}(\textcolor{black}{(N_h)_U(U_n)})}{\beta + \lambda_{\min}(A_h)}, & \text{if } \beta \geq \frac{\lambda_{\min}(\textcolor{black}{(N_h)_U(U_n)}) + \lambda_{\max}(\textcolor{black}{(N_h)_U(U_n)})}{2}.
\end{cases}
\]
  We compute the derivative \( g'(\beta) \) to find the minimum of \( g(\beta) \):

\[
g'(\beta) =
\begin{cases}
\displaystyle -\frac{\lambda_{\max}(\textcolor{black}{(N_h)_U(U_n)}) + \lambda_{\min}(A_h)}{(\beta + \lambda_{\min}(A_h))^2} < 0, & \text{if } -\lambda_{\min}(A_h) < \beta < \frac{\lambda_{\min}(\textcolor{black}{(N_h)_U(U_n)}) +\lambda_{\max}(\textcolor{black}{(N_h)_U(U_n)})}{2}, \\
\displaystyle \frac{\lambda_{\min}(\textcolor{black}{(N_h)_U(U_n)}) + \lambda_{\min}(A_h)}{(\beta + \lambda_{\min}(A_h))^2} > 0, & \text{if } \beta \geq \frac{\lambda_{\min}(\textcolor{black}{(N_h)_U(U_n)}) + \lambda_{\max}(\textcolor{black}{(N_h)_U(U_n)})}{2}.
\end{cases}
\]

From the sign of \( g'(\beta) \), we observe that \( g(\beta) \) decreases for \( \beta < \frac{\lambda_{\min}(\textcolor{black}{(N_h)_U(U_n)})+ \lambda_{\max}(\textcolor{black}{(N_h)_U(U_n)})}{2} \) and increases for \( \beta > \frac{\lambda_{\min}(\textcolor{black}{(N_h)_U(U_n)}) + \lambda_{\max}(\textcolor{black}{(N_h)_U(U_n)})}{2} \). Therefore, the minimum of \( g(\beta) \) occurs at:
\[
\beta = \frac{\lambda_{\min}(\textcolor{black}{(N_h)_U(U_n)})+ \lambda_{\max}(\textcolor{black}{(N_h)_U(U_n)})}{2}.
\]
\end{remark}
\begin{remark}
    \textcolor{black}{$N_h(\vec{x}, U_h)$ 
is sufficiently smooth to apply the Taylor expansion in \eqref{eq:F_u_n}, as \( f(x, \cdot) \in C^2(\mathbb{R}) \) for our model problem. Furthermore, from Remark \ref{Remark_integral}, each element of the vector  
$N_h(\vec{x}, U_h) = \operatorname{vec}([f(x_i, y_j, (U_h)_{i,j})])$
depends only on \( (x_i, y_j) \) and \( (U_h)_{i,j} \). Consequently, the matrix \( (N_h)_U \) exists and is diagonal, making the computation of its eigenvalues highly efficient.}
\end{remark}

\begin{remark}\label{remark}
We compute  the optimal upper bound given by the product \( \| A_h + \beta I \|^{-1} \| \beta I - (N_h)_U \| \), noting that it satisfies
\[
\| [A_h + \beta I]^{-1} (\beta I - (N_h)_U) \| \leq \| A_h + \beta I \|^{-1} \| \beta I - (N_h)_U \|.
\]
This approach offers a practical criterion for selecting
\( \beta \) to achieve the optimal upper bound. However, the value of \( \beta \) that minimizes \( \| [A_h + \beta I]^{-1} (\beta I - (N_h)_U) \| \) may differ from the proposed $\beta$ value.
\end{remark}

\section{Numerical Results}\label{examples}
In this section, all computations were performed using \textsc{Matlab} R2023b. For one-dimensional computations, we could not fully exploit the computational advantages of Graphics Processing Units (GPUs); therefore, all computation times are based on Central Processing Unit (CPU) computations using an AMD EPYC 7763 64-Core Processor. For two-dimensional and three-dimensional computations, we utilized GPUs, and all computation times are based on computations using an NVIDIA A100 80GB GPU.
\subsection{One-Dimensional Examples}\label{subsec:1D-examples}

In this sub-section, we present numerical examples in one dimension to illustrate the effectiveness of the quasi-Newton Method. In order to check the local convergence, we employ the initial guess as $U^*(x)+\varepsilon p(x)$,
where the perturbations are polynomial functions defined as:
\begin{equation*}
    p(x) = \frac{\sum_{i=0}^4 a_i x^i}{\left\| \sum_{i=0}^4 a_i x^i \right\|_\infty}, \quad x \in (0,1),
\end{equation*}
with the coefficients \( a_i \) are randomly selected from the interval \( (-1, 1) \), where \( \epsilon \) is the perturbation scale. The polynomial is normalized to ensure \( \| p(x) \|_\infty = 1 \). 
We conducted experiments using 100 random polynomials, computing the mean and standard deviation of both the iteration counts and computational times required for the convergence test.
More specifically, we compare our method with classical Newton's method when degree of freedom is $1000$. 

\subsubsection{Example 1}

First, we consider the following boundary value problem:
\begin{equation}\label{ex1}
\begin{cases}
    -u'' + u^3 = \pi^2 \sin(\pi x) + \sin^3(\pi x), & \text{in } (0,1), \\
    u(0) = u(1) = 0,
\end{cases}
\end{equation}
which has an analytic solution
$
    u(x) = \sin(\pi x).$
This example serves as a benchmark to validate the accuracy and efficiency of the proposed quasi-Newton method. We compared it with the classical Newton's method, which converges faster in this instance due to the absence of a tensor product implementation in 1D. In terms of convergence testing, the quasi-Newton method demonstrates comparable performance to Newton's method across various perturbation scales. Additionally, the proposed $\beta$ is shown to be optimal when compared to the quasi-Newton method with a fixed 
$\beta$.
The results are summarized in Table~\ref{table:results1}.

\begin{table}[!ht]
\centering
\begin{tabular}{|c|c|c|c|c|c|}
\hline
\multirow{2}{*}{Method} & \multirow{2}{*}{Metric} & \multicolumn{4}{c|}{Perturbation Scale} \\
\cline{3-6}
 & & 0.1 & 0.2 & 0.5 & 1.0  \\ \hline\multirow{2}{*}{Newtons method} & Iterations & 3.8 $\pm$ 0.4 & 4.0 $\pm$ 0.0 & 4.3 $\pm$ 0.5 & 4.7 $\pm$ 0.5\\ \cline{2-6} 
 & Time (s) & 0.02 $\pm$ 0.02 & 0.02 $\pm$ 0.01 & 0.02 $\pm$ 0.00 & 0.02 $\pm$ 0.00\\ \hline 
\multirow{2}{*}{$\beta=\frac{min+max}{2}$} & Iterations & 7.6 $\pm$ 0.5 & 7.8 $\pm$ 0.4 & 8.1 $\pm$ 0.5 & 8.3 $\pm$ 0.5\\ \cline{2-6} 
 & Time (s) & 0.07 $\pm$ 0.01 & 0.07 $\pm$ 0.01 & 0.08 $\pm$ 0.01 & 0.07 $\pm$ 0.01\\ \hline 
\multirow{2}{*}{$\beta={min+max}$} & Iterations & 7.7 $\pm$ 0.5 & 7.9 $\pm$ 0.3 & 8.2 $\pm$ 0.6 & 8.9 $\pm$ 0.6\\ \cline{2-6} 
 & Time (s) & 0.07 $\pm$ 0.01 & 0.07 $\pm$ 0.00 & 0.07 $\pm$ 0.01 & 0.08 $\pm$ 0.01\\ \hline 
\multirow{2}{*}{$\beta=max$} & Iterations & 7.8 $\pm$ 0.4 & 7.9 $\pm$ 0.3 & 8.2 $\pm$ 0.5 & 8.8 $\pm$ 0.5\\ \cline{2-6} 
 & Time (s) & 0.07 $\pm$ 0.01 & 0.07 $\pm$ 0.01 & 0.08 $\pm$ 0.01 & 0.07 $\pm$ 0.01\\ \hline 
\multirow{2}{*}{$\beta=min$} & Iterations & 12.4 $\pm$ 0.8 & 12.6 $\pm$ 0.7 & 13.4 $\pm$ 0.7 & 13.9 $\pm$ 0.8\\ \cline{2-6} 
 & Time (s) & 0.12 $\pm$ 0.01 & 0.11 $\pm$ 0.01 & 0.13 $\pm$ 0.01 & 0.12 $\pm$ 0.01\\ \hline 
\multirow{2}{*}{$\beta=0.1$} & Iterations & 11.9 $\pm$ 0.5 & 12.4 $\pm$ 0.7 & 12.8 $\pm$ 0.8 & 13.3 $\pm$ 0.5\\ \cline{2-6} 
 & Time (s) & 0.11 $\pm$ 0.01 & 0.12 $\pm$ 0.02 & 0.12 $\pm$ 0.01 & 0.11 $\pm$ 0.01\\ \hline 
\multirow{2}{*}{$\beta=0.5$} & Iterations & 10.6 $\pm$ 0.6 & 10.8 $\pm$ 0.7 & 11.4 $\pm$ 0.6 & 11.6 $\pm$ 0.7\\ \cline{2-6} 
 & Time (s) & 0.10 $\pm$ 0.01 & 0.09 $\pm$ 0.02 & 0.11 $\pm$ 0.01 & 0.09 $\pm$ 0.01\\ \hline 
\multirow{2}{*}{$\beta=1$} & Iterations & 8.9 $\pm$ 0.3 & 9.3 $\pm$ 0.5 & 9.5 $\pm$ 0.5 & 9.8 $\pm$ 0.8\\ \cline{2-6} 
 & Time (s) & 0.10 $\pm$ 0.01 & 0.08 $\pm$ 0.01 & 0.09 $\pm$ 0.01 & 0.19 $\pm$ 0.10\\ \hline 
\multirow{2}{*}{$\beta=5$} & Iterations & 10.8 $\pm$ 0.4 & 11.3 $\pm$ 0.6 & 11.8 $\pm$ 0.7 & 12.2 $\pm$ 0.7\\ \cline{2-6} 
 & Time (s) & 0.10 $\pm$ 0.01 & 0.09 $\pm$ 0.01 & 0.11 $\pm$ 0.01 & 0.12 $\pm$ 0.04\\ \hline 
\multirow{2}{*}{$\beta=10$} & Iterations & 18.3 $\pm$ 0.8 & 18.8 $\pm$ 1.0 & 19.8 $\pm$ 0.9 & 20.4 $\pm$ 1.1\\ \cline{2-6} 
 & Time (s) & 0.16 $\pm$ 0.02 & 0.21 $\pm$ 0.04 & 0.18 $\pm$ 0.01 & 0.19 $\pm$ 0.01\\ \hline 
\end{tabular}

\caption{Convergence comparisons between the quasi-Newton method and Newton's method for different values of 
$\beta$ and perturbation scales on the 1D Equation \eqref{ex1} are presented as means $\pm$ standard deviations across 100 random initializations.}
\label{table:results1}
\end{table}

\subsubsection{Example 2}
\begin{figure}[ht!]
\centering
\includegraphics[width=14cm,height=6cm,keepaspectratio]{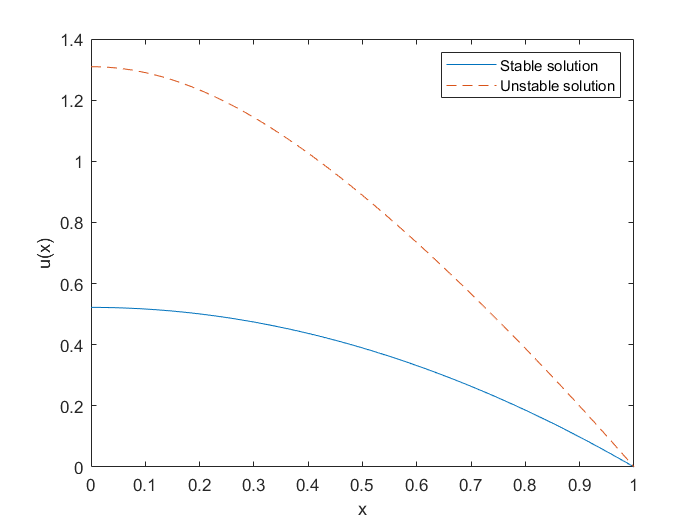}
   \caption{Numerical solutions of Eq. (\ref{ex2}) with $N=1000$ grid points. {\color{black}The unstable solution is plotted with dashed lines, while the stable solution is represented with solid lines.}}
   \label{1d_2}
\end{figure}

Next, we consider the following boundary value problem, which has two solutions. Here, we focus only on the stable solution, represented by solid lines in Fig. \ref{1d_2}. A comparison of the proposed quasi-Newton method shows comparable performance to Newton's method across various perturbation scales, as summarized in Table~\ref{table:results2}.
\begin{equation}\label{ex2}
\begin{cases} -u'' - (1+u^4) = 0, \quad \text{in } \quad (0,1), \\u'(0)=u(1)=0.\end{cases}
\end{equation}

\begin{table}[!ht]
\centering
\begin{tabular}{|c|c|c|c|c|c|}
\hline
\multirow{2}{*}{Method} & \multirow{2}{*}{Metric} & \multicolumn{4}{c|}{Perturbation Scale} \\
\cline{3-6}
 & & 0.1 & 0.2 & 0.3 & 0.4  \\ \hline\multirow{2}{*}{Newtons method} & Iterations & 3.9 $\pm$ 0.3 & 4.3 $\pm$ 0.5 & 4.3 $\pm$ 0.5 & 4.6 $\pm$ 1.0\\ \cline{2-6} 
 & Time (s) & 0.03 $\pm$ 0.01 & 0.02 $\pm$ 0.01 & 0.03 $\pm$ 0.01 & 0.03 $\pm$ 0.01\\ \hline 
\multirow{2}{*}{$\beta=\frac{min+max}{2}$} & Iterations & 6.6 $\pm$ 0.5 & 6.7 $\pm$ 0.5 & 6.9 $\pm$ 0.5 & 7.2 $\pm$ 0.6\\ \cline{2-6} 
 & Time (s) & 0.12 $\pm$ 0.02 & 0.11 $\pm$ 0.02 & 0.11 $\pm$ 0.02 & 0.11 $\pm$ 0.02\\ \hline 
\multirow{2}{*}{$\beta=max$} & Iterations & 9.3 $\pm$ 0.7 & 9.8 $\pm$ 0.6 & 10.1 $\pm$ 0.6 & 10.3 $\pm$ 0.8\\ \cline{2-6} 
 & Time (s) & 0.15 $\pm$ 0.03 & 0.15 $\pm$ 0.03 & 0.15 $\pm$ 0.03 & 0.15 $\pm$ 0.03\\ \hline 
\multirow{2}{*}{$\beta=0.1$} & Iterations & 10.2 $\pm$ 0.8 & 10.7 $\pm$ 0.6 & 10.9 $\pm$ 0.7 & 11.3 $\pm$ 0.7\\ \cline{2-6} 
 & Time (s) & 0.14 $\pm$ 0.03 & 0.15 $\pm$ 0.03 & 0.16 $\pm$ 0.03 & 0.16 $\pm$ 0.04\\ \hline 
\multirow{2}{*}{$\beta=0.5$} & Iterations & 13.3 $\pm$ 0.9 & 13.8 $\pm$ 1.0 & 14.5 $\pm$ 0.9 & 14.6 $\pm$ 1.0\\ \cline{2-6} 
 & Time (s) & 0.18 $\pm$ 0.03 & 0.20 $\pm$ 0.04 & 0.21 $\pm$ 0.04 & 0.20 $\pm$ 0.04\\ \hline 

\end{tabular}
\caption{
Convergence comparisons between the quasi-Newton method and Newton's method for different values of 
$\beta$ and perturbation scales on the 1D Equation \eqref{ex2} are presented as means $\pm$ standard deviations across 100 random initializations.}
\label{table:results2}
\end{table}

\textcolor{black}{In this example, since \( (N_h)_U \) is not SPD, we cannot guarantee that the upper bound for  
$
\| [A_h + \beta I]^{-1} \| \| (\beta I - (N_h)_U) \|$
is less than 1, as this is not theoretically provable without the SPD assumption. However, the proposed method still performs well empirically, as demonstrated by the numerical experiments. Our results indicate that the method continues to yield accurate and stable solutions even in the absence of the SPD condition.  
Additionally, it is important to note that the actual value of  
$
\| [A_h + \beta I]^{-1} (\beta I - (N_h)_U) \|$
could be smaller than the upper bound derived, as mentioned in Remark~\ref{remark}. The upper bound represents the worst-case scenario, whereas in practice,  
$
\| [A_h + \beta I]^{-1} (\beta I - (N_h)_U) \|$
may be significantly smaller than expected.}

\subsubsection{Example 3}
\begin{figure}[ht!]
\centering
\includegraphics[width=14cm,height=6cm,keepaspectratio]{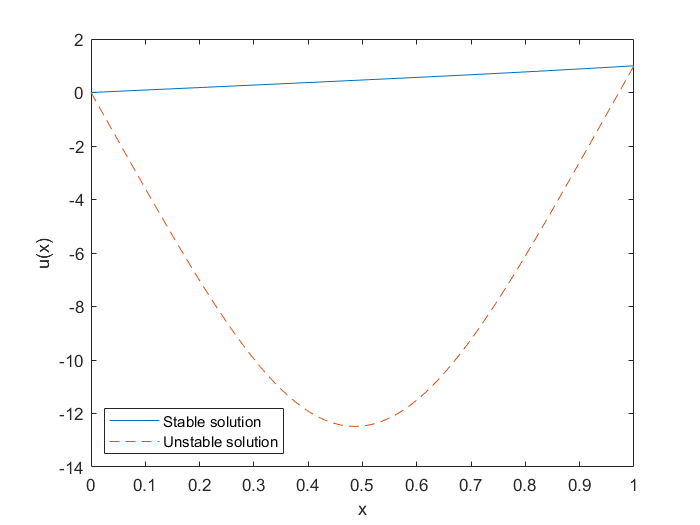}
   \caption{Numerical solutions of Eq. (\ref{example3}) with $N=1000$ grid points. {\color{black}The unstable solution is plotted with dashed lines, while the stable solution is represented with solid lines.}}
   \label{1d_3}
\end{figure}

Next, we consider the following boundary value problem, which also has multiple solutions. In this analysis, we focus only on the stable solution, represented by the solid lines in Fig.~\ref{1d_3}. The comparison results are presented in Table~\ref{table:results3}.
\begin{equation}\label{example3}
\begin{cases} - u'' + u^2 = 0, \quad \text{ in } \quad (0,1), \\u(0)=0\hbox{~and~}u(1)=1,\end{cases}
\end{equation}

\begin{table}[!ht]
\centering
\begin{tabular}{|c|c|c|c|c|c|}
\hline
\multirow{2}{*}{Method} & \multirow{2}{*}{Metric} & \multicolumn{4}{c|}{Perturbation Scale} \\
\cline{3-6}
 & & 0.1 & 0.2 & 0.5 & 1.0  \\ \hline\multirow{2}{*}{Newtons method} & Iterations & 3.5 $\pm$ 0.5 & 4.0 $\pm$ 0.2 & 4.0 $\pm$ 0.0 & 4.3 $\pm$ 0.5\\ \cline{2-6} 
 & Time (s) & 0.02 $\pm$ 0.00 & 0.01 $\pm$ 0.00 & 0.01 $\pm$ 0.00 & 0.01 $\pm$ 0.00\\ \hline 
\multirow{2}{*}{$\beta=\frac{min+max}{2}$} & Iterations & 6.0 $\pm$ 0.0 & 6.0 $\pm$ 0.0 & 6.0 $\pm$ 0.1 & 6.4 $\pm$ 0.5\\ \cline{2-6} 
 & Time (s) & 0.05 $\pm$ 0.01 & 0.05 $\pm$ 0.00 & 0.06 $\pm$ 0.01 & 0.06 $\pm$ 0.01\\ \hline 
\multirow{2}{*}{$\beta={min+max}$} & Iterations & 7.9 $\pm$ 0.3 & 8.4 $\pm$ 0.6 & 8.9 $\pm$ 0.4 & 8.9 $\pm$ 0.5\\ \cline{2-6} 
 & Time (s) & 0.07 $\pm$ 0.01 & 0.08 $\pm$ 0.01 & 0.09 $\pm$ 0.02 & 0.08 $\pm$ 0.01\\ \hline 
\multirow{2}{*}{$\beta=max$} & Iterations & 8.0 $\pm$ 0.3 & 8.6 $\pm$ 0.5 & 8.9 $\pm$ 0.4 & 9.2 $\pm$ 0.6\\ \cline{2-6} 
 & Time (s) & 0.07 $\pm$ 0.00 & 0.08 $\pm$ 0.02 & 0.10 $\pm$ 0.01 & 0.08 $\pm$ 0.01\\ \hline 
\multirow{2}{*}{$\beta=min$} & Iterations & 8.4 $\pm$ 0.6 & 8.8 $\pm$ 0.4 & 9.1 $\pm$ 0.4 & 9.6 $\pm$ 0.6\\ \cline{2-6} 
 & Time (s) & 0.07 $\pm$ 0.01 & 0.07 $\pm$ 0.01 & 0.10 $\pm$ 0.01 & 0.08 $\pm$ 0.01\\ \hline 
\multirow{2}{*}{$\beta=0.1$} & Iterations & 8.0 $\pm$ 0.2 & 8.3 $\pm$ 0.6 & 8.8 $\pm$ 0.5 & 8.7 $\pm$ 0.6\\ \cline{2-6} 
 & Time (s) & 0.07 $\pm$ 0.00 & 0.07 $\pm$ 0.01 & 0.09 $\pm$ 0.01 & 0.08 $\pm$ 0.01\\ \hline 
\multirow{2}{*}{$\beta=0.5$} & Iterations & 6.9 $\pm$ 0.3 & 6.9 $\pm$ 0.2 & 7.3 $\pm$ 0.5 & 7.7 $\pm$ 0.5\\ \cline{2-6} 
 & Time (s) & 0.06 $\pm$ 0.00 & 0.06 $\pm$ 0.00 & 0.07 $\pm$ 0.01 & 0.07 $\pm$ 0.01\\ \hline 
\multirow{2}{*}{$\beta=1$} & Iterations & 6.0 $\pm$ 0.0 & 6.0 $\pm$ 0.0 & 6.2 $\pm$ 0.4 & 6.5 $\pm$ 0.5\\ \cline{2-6} 
 & Time (s) & 0.06 $\pm$ 0.01 & 0.05 $\pm$ 0.00 & 0.06 $\pm$ 0.00 & 0.06 $\pm$ 0.00\\ \hline 
\multirow{2}{*}{$\beta=5$} & Iterations & 13.5 $\pm$ 0.9 & 14.1 $\pm$ 0.7 & 14.6 $\pm$ 0.9 & 15.4 $\pm$ 0.8\\ \cline{2-6} 
 & Time (s) & 0.14 $\pm$ 0.01 & 0.14 $\pm$ 0.02 & 0.14 $\pm$ 0.01 & 0.13 $\pm$ 0.01\\ \hline 
\multirow{2}{*}{$\beta=10$} & Iterations & 21.5 $\pm$ 0.9 & 22.1 $\pm$ 1.2 & 23.3 $\pm$ 1.0 & 24.2 $\pm$ 1.4\\ \cline{2-6} 
 & Time (s) & 0.20 $\pm$ 0.03 & 0.19 $\pm$ 0.02 & 0.22 $\pm$ 0.01 & 0.20 $\pm$ 0.01\\ \hline 
\end{tabular}

\caption{
Convergence comparisons between the quasi-Newton method and Newton's method for different values of 
$\beta$ and perturbation scales on the 1D Equation \eqref{example3} are presented as means $\pm$ standard deviations across 100 random initializations.}
\label{table:results3}
\end{table}

\subsection{Two dimensional Examples} 
Similarly to the 1D case, the perturbations in 2D are defined as polynomial functions:
\begin{equation*}
    p(x, y) = \frac{\sum_{i=0}^4 \sum_{j=0}^4 a_{ij} x^i y^j}{\left\| \sum_{i=0}^4 \sum_{j=0}^4 a_{ij} x^i y^j \right\|_\infty}, \quad (x, y) \in \Omega,
\end{equation*}
where the coefficients \( a_{ij} \) are randomly selected from the interval \( (-1, 1) \), with \( \epsilon \) as the perturbation scale. The normalization ensures that \( \| p(x, y) \|_\infty = 1 \). Our initial functions are generated by adding a random perturbation \(\epsilon p(x,y)\) to the numerical solution, allowing us to observe convergence behavior starting from perturbed initial conditions.
We also performed numerical experiments with 100 randomly generated polynomials, calculating the mean and standard deviation of the iteration counts and computational times required for the convergence test in each case.

\subsubsection{Example 1}
First, we consider the following Neumann boundary value problem:
\begin{equation}\label{ex2d1}
\begin{cases} -\Delta u(x,y)+u^3(x,y)+u(x,y)=-\Delta f(x,y)+f^3(x,y)+f(x,y)  \quad \text{in} \quad \Omega, \\ \frac{\partial u(x,y)}{\partial \textbf{n}}=0,\quad  \text{in} \quad \partial \Omega\end{cases} 
\end{equation}
where $\Omega=(-1,1) \times (-1,1)$ with $f(x,y)=\cos(\pi x)\cos (2\pi y)$ and \( \textbf{n} \) as the outward normal vector. We can easily find that the analytic solution is $u(x,y)=\cos(\pi x)\cos (2\pi y)$. The comparison results are presented in Table~\ref{table:results2d1} for {various} perturbation sizes and in Table~\ref{table:results2d_space} for different degree of freedom sizes. Table~\ref{table:results2d1} indicates that the Quasi-Newton method is {considerably more} efficient than Newton's method {on coarse meshes}, as it {leverages} the tensor product structure. {Additionally}, Table~\ref{table:results2d_space} demonstrates that the Quasi-Newton method {remains efficient even for finer mesh sizes}.

\begin{table}[!ht] 
 \centering 
 \begin{tabular}{|c|c|c|c|c|c|} 
 \hline 
\multirow{2}{*}{Method} & \multirow{2}{*}{Metric} & \multicolumn{4}{c|}{Perturbation Scale} \\ \cline{3-6}
 & & 0.1 & 0.2 & 0.5 & 1.0 \\ \hline\multirow{2}{*}{Newtons method} & Iterations & 4.0 $\pm$ 0.0 & 4.0 $\pm$ 0.0 & 4.6 $\pm$ 0.5 & 5.0 $\pm$ 0.2\\ \cline{2-6} 
 & GPU Time (s) & 0.4 $\pm$ 0.1 & 0.4 $\pm$ 0.1 & 0.5 $\pm$ 0.1 & 0.5 $\pm$ 0.1\\ \hline 
\multirow{2}{*}{$\beta=\frac{min+max}{2}$} & Iterations & 16.6 $\pm$ 0.8 & 17.4 $\pm$ 0.9 & 18.5 $\pm$ 0.9 & 19.6 $\pm$ 1.0\\ \cline{2-6} 
 & GPU Time (s) & 0.1 $\pm$ 0.0 & 0.1 $\pm$ 0.0 & 0.1 $\pm$ 0.0 & 0.1 $\pm$ 0.0\\ \hline 
\multirow{2}{*}{$\beta={min+max}$} & Iterations & 31.9 $\pm$ 1.5 & 33.1 $\pm$ 1.9 & 34.9 $\pm$ 2.1 & 36.8 $\pm$ 2.2\\ \cline{2-6} 
 & GPU Time (s) & 0.2 $\pm$ 0.1 & 0.2 $\pm$ 0.1 & 0.2 $\pm$ 0.1 & 0.2 $\pm$ 0.1\\ \hline 
\multirow{2}{*}{$\beta=max$} & Iterations & 31.4 $\pm$ 1.9 & 32.9 $\pm$ 1.9 & 35.0 $\pm$ 1.8 & 37.1 $\pm$ 2.1\\ \cline{2-6} 
 & GPU Time (s) & 0.2 $\pm$ 0.1 & 0.2 $\pm$ 0.1 & 0.2 $\pm$ 0.1 & 0.2 $\pm$ 0.1\\ \hline 
\multirow{2}{*}{$\beta=min$} & Iterations & 68.9 $\pm$ 3.4 & 70.8 $\pm$ 3.6 & 73.8 $\pm$ 4.8 & 77.2 $\pm$ 4.4\\ \cline{2-6} 
 & GPU Time (s) & 0.4 $\pm$ 0.1 & 0.4 $\pm$ 0.1 & 0.4 $\pm$ 0.1 & 0.4 $\pm$ 0.1\\ \hline 
\multirow{2}{*}{$\beta=3$} & Iterations & 21.5 $\pm$ 1.3 & 22.4 $\pm$ 1.0 & 23.3 $\pm$ 1.3 & 24.1 $\pm$ 1.2\\ \cline{2-6} 
 & GPU Time (s) & 0.1 $\pm$ 0.0 & 0.1 $\pm$ 0.0 & 0.1 $\pm$ 0.0 & 0.1 $\pm$ 0.0\\ \hline 
\multirow{2}{*}{$\beta=5$} & Iterations & 40.9 $\pm$ 2.6 & 42.4 $\pm$ 2.8 & 44.6 $\pm$ 2.2 & 46.3 $\pm$ 2.7\\ \cline{2-6} 
 & GPU Time (s) & 0.3 $\pm$ 0.1 & 0.3 $\pm$ 0.1 & 0.2 $\pm$ 0.1 & 0.2 $\pm$ 0.1\\ \hline 
\multirow{2}{*}{$\beta=7$} & Iterations & 60.5 $\pm$ 3.2 & 62.8 $\pm$ 3.1 & 66.6 $\pm$ 2.9 & 68.7 $\pm$ 2.8\\ \cline{2-6} 
 & GPU Time (s) & 0.4 $\pm$ 0.1 & 0.4 $\pm$ 0.1 & 0.4 $\pm$ 0.1 & 0.4 $\pm$ 0.1\\ \hline 
\multirow{2}{*}{$\beta=9$} & Iterations & 80.2 $\pm$ 3.9 & 83.1 $\pm$ 4.6 & 86.4 $\pm$ 5.7 & 90.5 $\pm$ 5.6\\ \cline{2-6} 
 & GPU Time (s) & 0.5 $\pm$ 0.1 & 0.5 $\pm$ 0.1 & 0.5 $\pm$ 0.1 & 0.5 $\pm$ 0.2\\ \hline 
\multirow{2}{*}{$\beta=11$} & Iterations & 98.0 $\pm$ 6.6 & 102.4 $\pm$ 6.8 & 107.5 $\pm$ 6.1 & 111.5 $\pm$ 6.9\\ \cline{2-6} 
 & GPU Time (s) & 0.6 $\pm$ 0.2 & 0.7 $\pm$ 0.1 & 0.6 $\pm$ 0.2 & 0.6 $\pm$ 0.2\\ \hline 
\end{tabular}
\caption{
Convergence comparisons between the quasi-Newton method and Newton's method for different values of 
$\beta$ and perturbation scales on the 2D Equation \eqref{ex2d1} are presented as means $\pm$ standard deviations across 100 random initialization
with $50 \times 50$ degree of freedom.}
\label{table:results2d1}
\end{table}

\begin{table}[!ht] 
 \centering 
 \begin{tabular}{|c|c|c|c|c|c|} 
 \hline 
\multirow{2}{*}{Method} & \multirow{2}{*}{Metric} & \multicolumn{4}{c|}{Degree of Freedom} \\ \cline{3-6}
 & & $200^2$  & $400^2$  & $800^2$ & $1600^2$ \\ \hline\multirow{2}{*}{$\beta=\frac{min+max}{2}$} & Iterations & 12.5 $\pm$ 0.9 & 12.7 $\pm$ 1.0 & 12.8 $\pm$ 0.8 & 12.9 $\pm$ 0.8\\ \cline{2-6} 
 & GPU Time (s) & 0.03 $\pm$ 0.00 & 0.03 $\pm$ 0.00 & 0.04 $\pm$ 0.01 & 0.10 $\pm$ 0.01\\ \hline 
\multirow{2}{*}{$\beta={min+max}$} & Iterations & 23.3 $\pm$ 1.9 & 23.2 $\pm$ 1.8 & 23.1 $\pm$ 2.0 & 23.8 $\pm$ 1.9\\ \cline{2-6} 
 & GPU Time (s) & 0.05 $\pm$ 0.01 & 0.05 $\pm$ 0.00 & 0.07 $\pm$ 0.01 & 0.17 $\pm$ 0.02\\ \hline 
\multirow{2}{*}{$\beta=max$} & Iterations & 23.4 $\pm$ 1.8 & 23.6 $\pm$ 1.7 & 23.8 $\pm$ 1.7 & 24.1 $\pm$ 1.7\\ \cline{2-6} 
 & GPU Time (s) & 0.05 $\pm$ 0.00 & 0.05 $\pm$ 0.00 & 0.07 $\pm$ 0.01 & 0.17 $\pm$ 0.01\\ \hline 
\multirow{2}{*}{$\beta=min$} & Iterations & 50.9 $\pm$ 4.0 & 51.3 $\pm$ 3.9 & 51.8 $\pm$ 3.1 & 52.4 $\pm$ 4.4\\ \cline{2-6} 
 & GPU Time (s) & 0.10 $\pm$ 0.01 & 0.11 $\pm$ 0.01 & 0.13 $\pm$ 0.01 & 0.33 $\pm$ 0.04\\ \hline 
\multirow{2}{*}{$\beta=1$} & Iterations & 51.5 $\pm$ 2.9 & 51.0 $\pm$ 4.3 & 51.1 $\pm$ 3.5 & 52.4 $\pm$ 3.7\\ \cline{2-6} 
 & GPU Time (s) & 0.09 $\pm$ 0.01 & 0.10 $\pm$ 0.01 & 0.12 $\pm$ 0.02 & 0.29 $\pm$ 0.02\\ \hline 
\multirow{2}{*}{$\beta=3$} & Iterations & 16.4 $\pm$ 1.0 & 16.2 $\pm$ 1.1 & 16.0 $\pm$ 1.4 & 16.2 $\pm$ 1.4\\ \cline{2-6} 
 & GPU Time (s) & 0.03 $\pm$ 0.00 & 0.03 $\pm$ 0.00 & 0.05 $\pm$ 0.01 & 0.09 $\pm$ 0.01\\ \hline 
\multirow{2}{*}{$\beta=5$} & Iterations & 30.0 $\pm$ 2.4 & 29.8 $\pm$ 3.1 & 30.4 $\pm$ 2.1 & 30.7 $\pm$ 2.7\\ \cline{2-6} 
 & GPU Time (s) & 0.06 $\pm$ 0.01 & 0.06 $\pm$ 0.01 & 0.07 $\pm$ 0.01 & 0.18 $\pm$ 0.02\\ \hline 
\multirow{2}{*}{$\beta=7$} & Iterations & 44.7 $\pm$ 4.0 & 45.0 $\pm$ 2.9 & 44.2 $\pm$ 3.6 & 45.5 $\pm$ 3.2\\ \cline{2-6} 
 & GPU Time (s) & 0.08 $\pm$ 0.01 & 0.09 $\pm$ 0.01 & 0.10 $\pm$ 0.01 & 0.26 $\pm$ 0.02\\ \hline 
\multirow{2}{*}{$\beta=9$} & Iterations & 58.6 $\pm$ 4.3 & 58.6 $\pm$ 4.8 & 58.2 $\pm$ 4.8 & 59.6 $\pm$ 4.8\\ \cline{2-6} 
 & GPU Time (s) & 0.11 $\pm$ 0.01 & 0.11 $\pm$ 0.01 & 0.12 $\pm$ 0.01 & 0.34 $\pm$ 0.03\\ \hline 
\end{tabular}
\caption{
Convergence comparisons between the quasi-Newton method and Newton's method for different values of 
$\beta$ and mesh sizes on the 2D Equation \eqref{ex2d1} are presented as means $\pm$ standard deviations across 100 random initialization with a perturbation size $0.1$.}
\label{table:results2d_space}
\end{table}

\subsubsection{Example 2}
\begin{figure}
    \centering
    \includegraphics[width=0.5\linewidth]{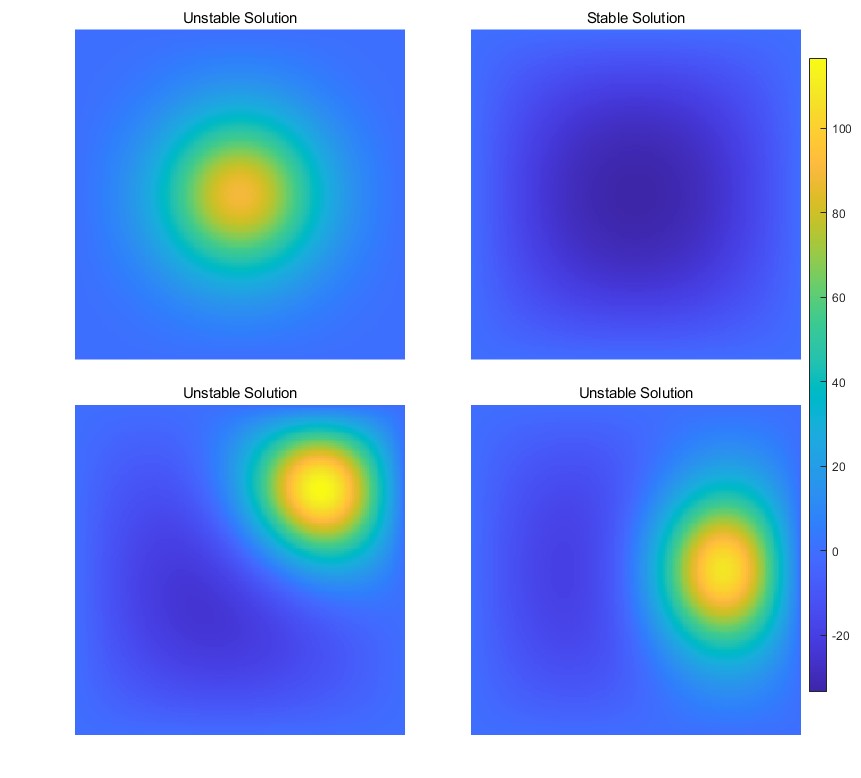}
    \caption{Four solutions of the Eq.~\eqref{2dex} on  $200 \times 200$ mesh.}
    \label{fig:2d_multi}
\end{figure}

Next, we considered the following  equation with multiple solutions:
\begin{equation}\label{2dex}
\begin{cases} 
-\Delta u(x,y) - u^2(x,y) = -s\sin(\pi x)\sin(\pi y), & \quad \text{in} \quad \Omega, \\
u(x,y) = 0, & \quad \text{on} \quad \partial \Omega 
\end{cases} 
\end{equation}
where $\Omega = (0,1) \times (0,1)$ and $s = 1600$. This example is adapted from \cite{breuer2003multiple}, there are four unique solutions up to rotation shown in Fig. \ref{fig:2d_multi}. 
We then focus exclusively on the stable solution. Table~\ref{table:results2D_2} presents the results obtained using the quasi-Newton method, demonstrating its efficiency across various perturbations and mesh sizes.

\begin{table}[!ht] 
 \centering 
 \begin{tabular}{|c|c|c|c|c|} 
 \hline 
\multirow{2}{*}{Metric} & \multicolumn{4}{c|}{Perturbation Scale} \\ \cline{2-5}
  & 0.1 & 0.2 & 0.5 & 1.0 \\ \hline Iterations & 24.0 $\pm$ 0.0 & 22.0 $\pm$ 0.0 & 25.0 $\pm$ 0.0 & 26.0 $\pm$ 0.0\\ \hline 
 GPU Time (s) & 0.04 $\pm$ 0.00 & 0.04 $\pm$ 0.01 & 0.04 $\pm$ 0.00 & 0.04 $\pm$ 0.00\\ \hline 
\multirow{2}{*}{Metric} & \multicolumn{4}{c|}{Degree of Freedom} \\ \cline{2-5}
 & $200^2$ & $400^2$ & $800^2$ & $1600^2$ \\ \hline Iterations & 24.0 $\pm$ 0.0 & 23.0 $\pm$ 0.0 & 21.0 $\pm$ 0.0 & 23.0 $\pm$ 0.0\\ \hline 
  GPU Time (s) & 0.04 $\pm$ 0.00 & 0.04 $\pm$ 0.00 & 0.03 $\pm$ 0.00 & 0.12 $\pm$ 0.00\\ \hline 
\end{tabular}
\caption{Convergence results of the quasi-Newton method for 2D Eq \eqref{2dex} with different perturbation size and mesh size, presented as means $\pm$ standard deviations across 100 random initialization. For different degree of freedom we choose perturbation size as $0.1$.}
\label{table:results2D_2}
\end{table}

\subsubsection{The Gray–Scott model in 2D}
\begin{figure}[!ht]
\centering
\includegraphics[width=16cm]{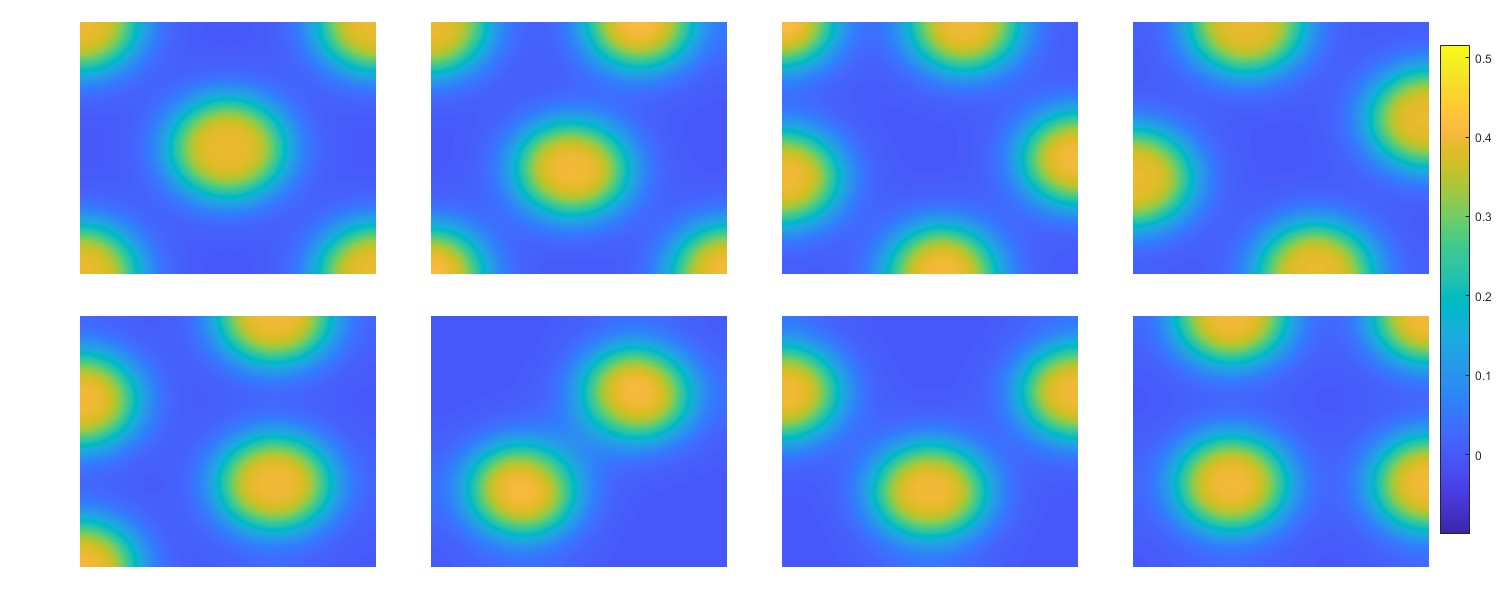}
 \caption{Eight solutions of the Gray-Scott model, showing only \( A(x,y) \) from the 2D system of Eq. \eqref{2dgrayeq} with a degree of freedom $1600 \times 1600$.}
 \label{2dgray}
\end{figure}

The last example for 2D is the steady-state Gray-Scott model, given by
\begin{equation}\label{2dgrayeq}
\begin{cases} -D_A\Delta A-SA^2+(\mu+\rho)A=0, \quad \text{in} \quad \Omega\\ -D_S\Delta S+SA^2-\rho(1-S)=0, \quad \text{in} \quad \Omega
\\ 
\displaystyle
\frac{\partial A}{\partial \textbf{n}}=\frac{\partial S}{\partial \textbf{n}}=0 \quad \text{in} \quad \partial \Omega

\end{cases}
\end{equation}
We consider the domain \( \Omega = (0,1) \times (0,1) \) with \( \textbf{n} \) as the outward normal vector \cite{hao2020spatial}. 
In this example, we set \( D_A = 2.5 \times 10^{-4} \), \( D_S = 5 \times 10^{-4} \), \( \rho = 0.04 \), and \( \mu = 0.065 \) to compute the multiple solutions.
%
We also set the step size to \( 0.1 \) in the quasi-Newton method to ensure convergence toward the non-trivial solution, starting from the initial guesses. The eight solutions computed by the quasi-Newton method are shown in Fig. \ref{2dgray}, and the total number of iterations and computation time are presented in Table \ref{table:results2D_g} for different mesh sizes.

\begin{table}[!ht] 
 \centering 
 \begin{tabular}{|c|c|c|c|c|} 
 \hline 
 \multirow{2}{*}{Metric} & \multicolumn{4}{c|}{Degree of Freedom} \\ \cline{2-5}
  & $200^2$ & $400^2$ & $800^2$& $1600^2$ \\ \hline Iterations  & 4505.9 & 4544.8 & 4550.0 & 4560.4\\ \hline
  GPU Time (s) &16.0 &  16.4 & 17.4 & 54.8\\ \hline
\end{tabular}
\caption{Convergence results for various mesh sizes in 2D Equation~\eqref{2dgray} with a perturbation size of \( 0.01 \). For each solution, the time to converge was computed individually; the total average time was then obtained by summing these times and dividing by 8. We report the total average time and iteration count required to converge to solutions.}
\label{table:results2D_g}
\end{table}
\textcolor{black}{
We have extended the proposed method to nonlinear differential equation systems, as discussed in \(\S \ref{ext_nonli_sys}\), and the results demonstrate promising performance. Despite the added complexity, the method consistently produces stable and accurate solutions, highlighting its capability and robustness in handling nonlinear dynamics effectively.
}

\subsection{Three dimensional Examples}
The perturbations are defined as polynomial functions:
\begin{equation*}
    p(x, y,z) = \frac{\sum_{i=0}^4 \sum_{j=0}^4\sum_{k=0}^4 a_{ijk} x^i y^j z^k}{\left\| \sum_{i=0}^4 \sum_{j=0}^4\sum_{k=0}^4 a_{ijk} x^i y^j z^k \right\|_\infty}, \quad (x, y, z) \in \Omega,
\end{equation*}
where the coefficients \( a_{ijk} \) are randomly selected from the interval \( (-1, 1) \), with \( \epsilon \) as the perturbation scale. The normalization ensures that \( \| p(x, y, z) \|_\infty = 1 \). Our initial functions are generated by adding a random perturbation \(\epsilon p(x,y,z)\) to the numerical solution, allowing us to observe convergence behavior starting from perturbed initial conditions.

We performed experiments with 100 randomly generated polynomials, calculating the mean and standard deviation of the iteration counts and computational times required for convergence in each case.

\subsubsection{Example 1}
First, we consider the following boundary value problem
\begin{equation}\label{ex3d1}
\begin{cases} -\triangle u+ a u +b u^3+V u = -\triangle f+ a f +b f^3+V f, \quad \text{in } \quad (-1,1) \times (-1,1) \times (-1,1), \\ \partial_nu=0  \quad \text{in } \quad \partial \Omega,\end{cases}
\end{equation}
where \( V \) is defined as
\[
V = 100\left( \sin^2\left( \frac{\pi x}{4} \right) + \sin^2\left( \frac{\pi y}{4} \right) + \sin^2\left( \frac{\pi z}{4} \right) \right) + x^2 + y^2 + z^2,
\] and $a=b=10$. With explicit solution as $u=f=\cos(\pi x)+\cos(\pi y)+\cos(\pi z)$. The comparison results are presented in Tables \ref{table:results3D} and \ref{table:results3D_mesh}, showing the optimal \(\beta\) for different perturbation sizes and mesh sizes.

\begin{table}[!ht] 
 \centering 
 \begin{tabular}{|c|c|c|c|c|c|} 
 \hline 
\multirow{2}{*}{Method} & \multirow{2}{*}{Metric} & \multicolumn{4}{c|}{Perturbation Scale} \\ \cline{3-6}
 & & 0.1 & 0.2 & 0.5 & 1.0 \\ \hline\multirow{2}{*}{$\beta=\frac{min+max}{2}$} & Iterations & 62.9 $\pm$ 1.6 & 65.0 $\pm$ 1.7 & 67.4 $\pm$ 1.6 & 69.7 $\pm$ 1.5\\ \cline{2-6} 
 & GPU Time (s) & 0.4 $\pm$ 0.0 & 0.4 $\pm$ 0.0 & 0.4 $\pm$ 0.0 & 0.5 $\pm$ 0.0\\ \hline 
\multirow{2}{*}{$\beta={min+max}$} & Iterations & 129.0 $\pm$ 2.8 & 133.1 $\pm$ 2.5 & 138.5 $\pm$ 2.7 & 142.2 $\pm$ 3.1\\ \cline{2-6} 
 & GPU Time (s) & 0.8 $\pm$ 0.0 & 0.8 $\pm$ 0.0 & 0.8 $\pm$ 0.0 & 0.9 $\pm$ 0.0\\ \hline 
\multirow{2}{*}{$\beta=max$} & Iterations & 116.8 $\pm$ 2.4 & 120.4 $\pm$ 3.0 & 125.4 $\pm$ 2.7 & 128.8 $\pm$ 2.4\\ \cline{2-6} 
 & GPU Time (s) & 0.7 $\pm$ 0.0 & 0.7 $\pm$ 0.0 & 0.7 $\pm$ 0.0 & 0.7 $\pm$ 0.0\\ \hline 
\multirow{2}{*}{$\beta=175$} & Iterations & 124.0 $\pm$ 2.4 & 127.6 $\pm$ 2.4 & 131.3 $\pm$ 2.5 & 134.2 $\pm$ 2.5\\ \cline{2-6} 
 & GPU Time (s) & 0.7 $\pm$ 0.0 & 0.7 $\pm$ 0.0 & 0.7 $\pm$ 0.0 & 0.7 $\pm$ 0.0\\ \hline 
\multirow{2}{*}{$\beta=200$} & Iterations & 63.0 $\pm$ 1.3 & 64.7 $\pm$ 1.5 & 66.2 $\pm$ 1.2 & 67.8 $\pm$ 1.5\\ \cline{2-6} 
 & GPU Time (s) & 0.4 $\pm$ 0.0 & 0.4 $\pm$ 0.0 & 0.4 $\pm$ 0.0 & 0.4 $\pm$ 0.0\\ \hline 
\multirow{2}{*}{$\beta=250$} & Iterations & 67.7 $\pm$ 1.5 & 69.9 $\pm$ 1.7 & 72.4 $\pm$ 1.7 & 74.5 $\pm$ 1.8\\ \cline{2-6} 
 & GPU Time (s) & 0.4 $\pm$ 0.0 & 0.4 $\pm$ 0.0 & 0.4 $\pm$ 0.0 & 0.4 $\pm$ 0.0\\ \hline 
\multirow{2}{*}{$\beta=300$} & Iterations & 82.1 $\pm$ 2.0 & 84.4 $\pm$ 1.9 & 87.8 $\pm$ 1.9 & 90.1 $\pm$ 1.8\\ \cline{2-6} 
 & GPU Time (s) & 0.5 $\pm$ 0.0 & 0.5 $\pm$ 0.0 & 0.5 $\pm$ 0.0 & 0.5 $\pm$ 0.0\\ \hline 
\multirow{2}{*}{$\beta=350$} & Iterations & 95.8 $\pm$ 2.1 & 98.7 $\pm$ 2.3 & 103.2 $\pm$ 1.8 & 105.8 $\pm$ 2.0\\ \cline{2-6} 
 & GPU Time (s) & 0.5 $\pm$ 0.0 & 0.5 $\pm$ 0.0 & 0.6 $\pm$ 0.0 & 0.6 $\pm$ 0.0\\ \hline 
\end{tabular}
\caption{Convergence results of different methods for 3D Eq. \eqref{ex3d1} with $200^3$ degree of freedom., presented as means $\pm$ standard deviations across 100 random initialization.}
\label{table:results3D}
\end{table}

\begin{table}[!ht] 
 \centering 
 \begin{tabular}{|c|c|c|c|c|} 
 \hline 
\multirow{2}{*}{Method} & \multirow{2}{*}{Metric} & \multicolumn{3}{c|}{Degree of Freedom} \\ \cline{3-5}
 & & $200^3$ & $400^3$ & $800^3$ \\ \hline\multirow{2}{*}{$\beta=\frac{min+max}{2}$} & Iterations & 62.9 $\pm$ 1.8 & 59.4 $\pm$ 1.6 & 55.5 $\pm$ 1.3\\ \cline{2-5} 
 & GPU Time (s) & 0.5 $\pm$ 0.0 & 3.5 $\pm$ 0.1 & 36.5 $\pm$ 0.9\\ \hline 
\multirow{2}{*}{$\beta={min+max}$} & Iterations & 128.8 $\pm$ 2.6 & 121.1 $\pm$ 2.3 & 113.1 $\pm$ 2.9\\ \cline{2-5} 
 & GPU Time (s) & 0.8 $\pm$ 0.0 & 6.6 $\pm$ 0.2 & 70.0 $\pm$ 1.7\\ \hline 
\multirow{2}{*}{$\beta=max$} & Iterations & 117.0 $\pm$ 2.4 & 109.8 $\pm$ 2.6 & 103.1 $\pm$ 2.5\\ \cline{2-5} 
 & GPU Time (s) & 0.7 $\pm$ 0.0 & 5.8 $\pm$ 0.2 & 62.4 $\pm$ 1.6\\ \hline 
\multirow{2}{*}{$\beta=175$} & Iterations & 124.3 $\pm$ 2.5 & 118.9 $\pm$ 2.3 & 113.0 $\pm$ 2.4\\ \cline{2-5} 
 & GPU Time (s) & 5.3 $\pm$ 0.0 & 47.9 $\pm$ 0.2 & 577.3 $\pm$ 3.8\\ \hline 
\multirow{2}{*}{$\beta=200$} & Iterations & 62.8 $\pm$ 1.3 & 60.2 $\pm$ 1.1 & 57.5 $\pm$ 1.2\\ \cline{2-5} 
 & GPU Time (s) & 0.4 $\pm$ 0.0 & 3.2 $\pm$ 0.1 & 35.9 $\pm$ 0.7\\ \hline 
\multirow{2}{*}{$\beta=250$} & Iterations & 67.7 $\pm$ 1.8 & 63.8 $\pm$ 1.5 & 59.7 $\pm$ 1.6\\ \cline{2-5} 
 & GPU Time (s) & 0.4 $\pm$ 0.0 & 3.4 $\pm$ 0.1 & 37.2 $\pm$ 0.9\\ \hline 
\multirow{2}{*}{$\beta=300$} & Iterations & 82.3 $\pm$ 1.8 & 77.5 $\pm$ 1.8 & 72.3 $\pm$ 1.7\\ \cline{2-5} 
 & GPU Time (s) & 0.5 $\pm$ 0.0 & 4.1 $\pm$ 0.1 & 44.4 $\pm$ 1.0\\ \hline 
\multirow{2}{*}{$\beta=350$} & Iterations & 96.2 $\pm$ 2.1 & 90.3 $\pm$ 2.1 & 84.6 $\pm$ 2.3\\ \cline{2-5} 
 & GPU Time (s) & 0.6 $\pm$ 0.0 & 4.7 $\pm$ 0.1 & 51.5 $\pm$ 1.3\\ \hline 
\end{tabular}
\caption{Convergence results for different methods with various mesh size in 3D Eq. \eqref{ex3d1} with $200^3$ degree of freedom, presented as means $\pm$ standard deviations across 100 random initialization. Here perturbation size is $0.1$.}
\label{table:results3D_mesh}
\end{table}

\subsubsection{Example 2}\label{3d_eigen_p}
\begin{figure}
    \centering
    \includegraphics[width=0.5\linewidth]{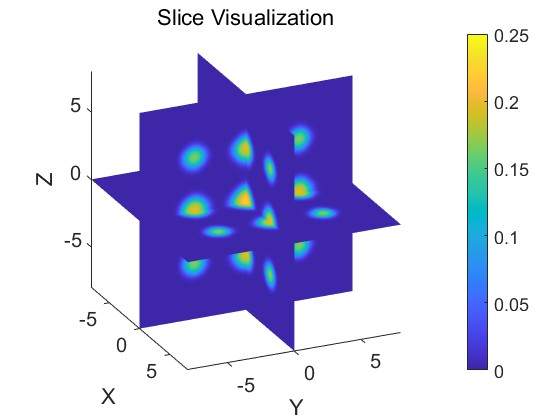}
    \caption{The 3D solution of Eq. \eqref{ex3d2}, with slices shown for the planes $x=0$, $y=0$ and $z=0$.}
    \label{3d_fig}
\end{figure}

Lastly, we consider the following Schrodinger equation:
\begin{equation}\label{ex3d2}
\begin{cases} 
    -\triangle u + \alpha |u|^2 u + V u = \lambda u, & \text{in } \quad (-8,8) \times (-8,8) \times (-8,8), \\ 
     u = 0, & \text{on } \quad \partial \Omega,
\end{cases}
\end{equation}
where \( V \) is defined as
\[
V = 100\left( \sin^2\left( \frac{\pi x}{4} \right) + \sin^2\left( \frac{\pi y}{4} \right) + \sin^2\left( \frac{\pi z}{4} \right) \right) + x^2 + y^2 + z^2,
\]
and \( \alpha = 1600 \).
This example is based on the model in \cite{chen2024fully}, where the value of \( \lambda \) is not provided explicitly. To compute \( \lambda \), we use the following formula:
\begin{equation*}
    \lambda = \frac{\int_{\Omega}(-\Delta u + f(u))u}{\int_\Omega (u)(u)},
\end{equation*}
and substituting this expression back into equation~\eqref{ex3d2}. After each Quasi-Newton method we normalized $u$. The detailed explanation is in \S \ref{explanation_eigp}. The results for this example by the proposed Quasi-Newton method are presented in Table~\ref{table:results3D_2} and Fig.~\ref{3d_fig} which show the efficiency with different perturbations and mesh sizes.

\begin{table}[!ht] 
 \centering 
 \begin{tabular}{|c|c|c|c|} 
 \hline 
\multirow{2}{*}{Metric} & \multicolumn{3}{c|}{Perturbation Scale} \\ \cline{2-4}
  & 0.1 & 0.2 & 0.3 \\ \hline Iterations & 761.4 $\pm$ 4.0 & 781.7 $\pm$ 5.5 & 773.7 $\pm$ 20.3 \\ \hline 
 GPU Time (s) & 4.6 $\pm$ 0.1 & 4.7 $\pm$ 0.1 & 16.0 $\pm$ 44.6  \\ \hline 
\multirow{2}{*}{Metric} & \multicolumn{3}{c|}{Degree of Freedom} \\ \cline{2-4}
 & $200^3$ & $400^3$ & $800^3$  \\ \hline Iterations & 761.4 $\pm$ 4.0 &704.0 $\pm$ 1.4 & 685.5 $\pm$ 2.1 \\ \hline 
 GPU Time (s) & 4.6 $\pm$ 0.1 & 37.4 $\pm$ 0.1 & 423.2 $\pm$ 1.2 \\ \hline 
\end{tabular}
\caption{Convergence results of the quasi-Newton method for 3D Eq. \eqref{ex3d2} with different perturbation and mesh sizes, presented as means $\pm$ standard deviations across 100 random initialization. For different mesh sizes, we choose perturbation size as $0.1$.}
\label{table:results3D_2}
\end{table}

\section{Conclusion}

In this paper, we have developed a quasi-Newton method for solving quasi-linear elliptic equations that optimally exploit GPU architectures for computational efficiency. By approximating the Jacobian matrix as a combination of a linear Laplacian and simplified nonlinear terms, our approach mitigates the high computational costs often associated with traditional Newton's methods, especially when dealing with large, sparse matrices from discretized PDEs. The method is highly adaptable for parallelized computations across tensor structures, making it ideal for tackling high-dimensional problems.

Our convergence analysis provides a theoretical foundation for the method, showing that it achieves local convergence to the exact solution \( U^\ast \) under appropriate choices for the regularization parameter \( \beta_n \). By setting \( \beta_n \) to balance the maximum and minimum eigenvalues of the nonlinear Jacobian, we ensure both stability and efficiency in each iteration. 

Practical applications demonstrate the robustness and scalability of our quasi-Newton method, especially in two- and three-dimensional spaces, and show its efficacy using the Spectral Element Method (SEM) and Finite Difference Method (FDM) on GPUs. This combination of theoretical and practical advantages makes our method a powerful tool for solving complex PDEs in physics, engineering, and related fields. Future work will explore extending this approach to other classes of nonlinear PDEs and refining the regularization strategy to adapt dynamically, further leveraging GPU capabilities for high-performance computing in scientific applications.

\renewcommand{\ackname}{Funding}

\begin{acknowledgements}
S.L, and W.H. is supported by both NIH via 1R35GM146894 and NSF DMS-2052685; 
 X.Z. is supported by NSF DMS-2208518.  
\end{acknowledgements}

\renewcommand{\ackname}{Data availability}

\begin{acknowledgements}
Data sharing is not applicable to this article as no datasets were generated or analyzed during
 the current study.
\end{acknowledgements}

\section*{Declarations}
\renewcommand{\ackname}{Conflict of interest}

\begin{acknowledgements}
The authors declare that they have no Conflict of interest.
\end{acknowledgements}

\bibliographystyle{plain}
\bibliography{references}

\end{document}